\documentclass[a4paper]{article}
\usepackage{amssymb}
\usepackage{amsfonts}
\usepackage{amsmath}

\input{tcilatex}

\begin{document}

\title{The characteristic classes and Weyl invariants of Spinor groups}
\author{Haibao Duan\thanks{%
The author's research is supported by NSFC 11131008; 11661131004} \\
Academy of Mathematics and Systems Science, \\
Chinese Academy of Sciences, Beijing 100190, China }
\maketitle

\begin{abstract}
Based on a pair of new cohomology operations on so called $\delta _{2}$%
--formal spaces, we construct the integral cohomologies of the classifying
spaces of the Lie groups $Spin(n)$ and $Spin^{c}(n)$.

As applications, we determine the ring of integral Weyl invariants of the
group $Spin(n)$, introduce the spin characteristic classes for the spin
vector bundles, and demonstrate their usages in spin geometry.

\begin{description}
\item[2010 Mathematical Subject Classification] 55R240, 55S20, 53C27

\item[Key words and phrases:] Spinor groups, classifying spaces,
characteristic classes, Weyl invariants

\item[Email:] dhb@math.ac.cn
\end{description}
\end{abstract}

\section{Introduction}

The spin group $Spin(n)$ is the universal cover of the special orthogonal
group $SO(n)$. The spin$^{c}$ group $Spin^{c}(n)$ is the central extension $%
Spin(n)\times _{\mathbb{Z}_{2}}U(1)$ of $SO(n)$ by the circle group $U(1)$.
In the first part of the paper we introduce a pair $F=(\alpha ,\gamma )$ of
secondary cohomology operations, which is applied to construct the integral
cohomology rings of the classifying spaces $B_{Spin^{c}(n)}$ and $%
B_{Spin(n)} $.

The $\func{mod}2$ cohomology of the space $B_{Spin(n)}$ has been determined
by Borel \cite{B1} for $n\leq 10$, and completed by Quillen in \cite{Q}.
Concerning the integral cohomology $H^{\ast }(B_{Spin(n)})$ partial
information are known \cite{Web,Web1}. In \cite{Th} Thomas calculated the
cohomology $H^{\ast }(B_{Spin(n)})$ in the stable range $n=\infty $, but his
result relies on two sequences $\{\Phi _{i}\},\{\Psi _{i}\}$ of
indeterminacies. Another inspiring approach is due to Benson and Wood. By
computing with the Weyl invariants a partial presentation of the ring $%
H^{\ast }(B_{Spin(n)})$ is formulated in \cite[Theorem 11.1]{BW}, where the
determination of explicit generators and relations is noted to be a rather
daunting task. For the difficulties that one encounters when computing with
the cohomologies of the classifying space $B_{G}$ of a Lie group $G$, we
refer to Feshbach \cite[Final remarks]{Fe}. In our approach the pair $%
F=(\alpha ,\gamma )$ of cohomology operators will make the structure of the
ring $H^{\ast }(B_{Spin(n)})$ appearing in a new light, see Remarks 8.6 and
9.5.

Knowing the integral cohomology of the classifying space $B_{G}$ of a Lie
group $G$ has direct consequences in geometry and invariant theory. In
particular, assuming that a minimal system $\{q_{1},\cdots ,q_{m}\}$ of
generators of the ring $H^{\ast }(B_{G})$ has been specified, one can
introduce the characteristic classes for a principle $G$ bundle $\xi $ over
a space $X$ by letting

\begin{quote}
$q_{r}(\xi ):=f_{\xi }^{\ast }(q_{r})\in H^{\ast }(X)$, $1\leq r\leq m$,
\end{quote}

\noindent where $f_{\xi }:X\rightarrow B_{G}$ is the classifying map of the
bundle $\xi $. One obtains also the basic Weyl invariants of the group $G$
by setting

\begin{quote}
$d_{r}:=B_{t}^{\ast }(q_{r})\in H^{\ast }(B_{T})$, $1\leq r\leq m$,
\end{quote}

\noindent where $T$ is a maximal torus on $G$, and the map $%
B_{t}:B_{T}\rightarrow B_{G}$ is induced by the inclusion $T$ $\subset G$.
For the classical groups $G=U(n),SO(n)$ and $Sp(n)$ these stories have been
well understood by the 1950's \cite{B,BH}. In the second part of the paper
we complete the projects for the spinor groups $Spin(n)$ and $Spin^{c}(n)$.

In mathematical physics the Postnikov tower anchored by the classifying
space $B_{SO(n)}$ is

\begin{quote}
$\cdots \rightarrow $ $B_{Fivebrane(n)}\rightarrow B_{String(n)}\rightarrow
B_{Spin(n)}\rightarrow B_{SO(n)}$,
\end{quote}

\noindent indicating that the calculation of ring $H^{\ast }(B_{Spin(n)})$
is a necessary step towards the integral cohomologies of the further spaces $%
B_{String(n)}$ and $B_{Fivebrane(n)}$ in the tower. In addition, the
integral cohomology $H^{\ast }(B_{Spin(n)})$ is essential to form the
cohomology of the Thom spectrum $\{M_{Spin(n)},n\geq 3\}$ \cite{ABP}, is
explicitly requested by the study of the complex cobordism of the
classifying space $B_{Spin(n)}$ \cite[Corollary 1.4]{T}, and by the
understanding of the total Chern class of the complex spin presentation of
the group $Spin(n)$ \cite{ABS,Q}.

\bigskip

\noindent \textbf{Remark 1.1.} If $3\leq n\leq 6$ there hold the following
group isomorphisms

\begin{quote}
$Spin(3)=SU(2)$, $Spin(4)=SU(2)\times SU(2)$,

$Spin(5)=Sp(2)$, $Spin(6)=SU(4)$,
\end{quote}

\noindent where $SU(k)$ is the special unitary group of rank $k$, $Sp(2)$ is
the symplectic group of rank $2$. Therefore, we may assume $n\geq 7$ to
exclude these cases.$\square $

\section{The main results}

For a topological space $X$ let $Sq^{k}$ be the Steenrod squares on the $%
\func{mod}2$ cohomology algebra $H^{\ast }(X;\mathbb{Z}_{2})$, and denote by 
$\delta _{m}$ the Bockstein homomorphism from the $\func{mod}m$ cohomology $%
H^{r}(X;\mathbb{Z}_{m})$ into the integral cohomology $H^{r+1}(X)$. For the
homomorphisms of coefficients groups

\begin{quote}
$\theta :\mathbb{Z}_{2}\rightarrow \mathbb{Z}_{4}$ by $\theta (1)=2$, and $%
\rho _{m}:\mathbb{Z}\rightarrow \mathbb{Z}_{m}$ by $\rho _{m}(1)=1$,
\end{quote}

\noindent the same notion are applied to denote their induced maps on the
cohomologies.

Let $\mathcal{B}:$ $H^{2r}(X;\mathbb{Z}_{2})\rightarrow H^{4r}(X;\mathbb{Z}%
_{4})$ be the Pontryagin square \cite{BT}. For an even degree class $u\in
H^{2r}(X;\mathbb{Z}_{2})$ there holds the following universal relation

\begin{enumerate}
\item[(2.1)] $\delta _{2}(u\cup u)=2\delta _{4}\mathcal{B}(u)$ on $%
H^{4r+1}(X)$ (see (3.1)).
\end{enumerate}

\noindent \textbf{Definition 2.1.} The space $X$ is called $\delta _{2}$%
\textsl{--formal} if $\delta _{2}(u\cup u)=0$ for every $u\in H^{2r}(X;%
\mathbb{Z}_{2})$, $r\geq 1$.$\square $

\bigskip

It follows from (2.1) that, if $X$\ is a space whose integral cohomologies $%
H^{\ast }(X)$ in degrees $4r+1$ have no torsion element of order $4$, then $%
X $\ is $\delta _{2}$--formal. In particular, all the $1$--connected Lie
groups, the classifying spaces $B_{SO(n)}$, $B_{Spin(n)}$ and $%
B_{Spin^{c}(n)}$ of our concern, as well as the Thom spectrum $%
\{M_{Spin(n)},n\geq 3\}$, are all examples of $\delta _{2}$--formal spaces.

To introduce the promised operations observe that the Bockstein operator $%
Sq^{1}=\rho _{2}\circ \delta _{2}$ defines the following decomposition on
the $\mathbb{Z}_{2}$ space $H^{\ast }(X;\mathbb{Z}_{2})$

\begin{quote}
$H^{\ast }(X;\mathbb{Z}_{2})=\ker Sq^{1}\oplus S_{2}^{\ast }(X)$ with $%
S_{2}^{\ast }(X)=H^{\ast }(X;\mathbb{Z}_{2})/\ker Sq^{1}$.
\end{quote}

\noindent \textbf{Theorem A.} \textsl{For any }$\delta _{2}$\textsl{--formal
space }$X$\textsl{\ there exists a unique pair of cohomological operations}

\begin{quote}
$F:H^{2r}(X;\mathbb{Z}_{2})\rightarrow H^{4r}(X;\mathbb{Z}_{4})\times
S^{4r}(X;\mathbb{Z}_{2})$\textsl{,}
\end{quote}

\noindent \textsl{written }$F(u)=(\alpha (u),\gamma (u))$\textsl{,} $u\in
H^{2r}(X;\mathbb{Z}_{2})$\textsl{, that is characterized by the following
three properties:}

\begin{quote}
\textsl{i) }$\alpha (u)\in \func{Im}\rho _{4}$\textsl{;}

\textsl{ii)} $\mathcal{B}(u)=\alpha (u)+\theta (\gamma (u))$\textsl{;}

\textsl{iii)} $Sq^{1}(\gamma (u))=Sq^{2r}Sq^{1}(u)+u\cup Sq^{1}(u)$\textsl{.}
\end{quote}

The uniqueness assertion of Theorem A implies that the properties i), ii)
and iii) can be taken as \textsl{an axiomatic definition} of the pair $%
F=(\alpha ,\gamma )$ of operators. In particular, since $Sq^{1}$ injects on $%
S_{2}^{\ast }(X)$ while $\gamma (u)\in S_{2}^{\ast }(X)$, the operator $%
\gamma $ is determined uniquely by the relation iii). Thus, applying $Sq^{1}$
to both sides verifies at once the following equalities useful to evaluate $%
\gamma $.

\bigskip

\noindent \textbf{Corollary 2.2.} \textsl{For any }$\delta _{2}$\textsl{%
--formal space }$X$\textsl{\ and }$u_{1}$\textsl{,}$u_{2}\in H^{\ast }(X;%
\mathbb{Z}_{2})$ \textsl{with} $\deg u_{i}=2r_{i}$\textsl{,} \textsl{one has}

\begin{quote}
\textsl{i)} $\gamma (u_{1}+u_{2})\equiv \gamma (u_{1})+\gamma
(u_{2})+u_{1}\cup u_{2}\func{mod}\ker Sq^{1}$

\textsl{ii)} $\gamma (u_{1}\cup u_{2})\equiv u_{1}^{2}{}\cup \gamma
(u_{2})+u_{2}^{2}{}\cup \gamma (u_{1})+u_{1}\cup Sq^{1}u_{1}\cup
Sq^{2r_{2}-1}u_{2}$

$\qquad +u_{2}\cup Sq^{1}u_{2}\cup Sq^{2r_{1}-1}u_{1}\func{mod}\ker Sq^{1}$.$%
\square $
\end{quote}

The operator $\gamma $ can be iterated to yield the following notion.

\bigskip

\noindent \textbf{Definition 2.3. }For an even degree cohomology class $u\in
H^{2r}(X;\mathbb{Z}_{2})$ of a $\delta _{2}$--formal space $X$, the sequence 
$\left\{ u^{(0)},u^{(1)},u^{(2)},\cdots \right\} $ on $H^{\ast }(X;\mathbb{Z}%
_{2})$ defined recursively by $u^{(0)}=u$, $u^{(k+1)}=\gamma (u^{(k)})$ is
called \textsl{the derived sequence} \textsl{of} $u$.$\square $

\bigskip

\noindent \textbf{Example 2.4.} Recall that $H^{\ast }(B_{SO(n)};\mathbb{Z}%
_{2})=\mathbb{Z}_{2}[w_{2},\cdots ,w_{n}]$, where $w_{i}$ is the $i^{th}$
Stiefel--Whitney class of the canonical real $n$--bundle on $B_{SO(n)}$. For 
$u=w_{2r}$ solving the equation iii) of Theorem A using
coefficient--comparison yields

\begin{enumerate}
\item[(2.2)] $\gamma (w_{2r})=w_{4r}+w_{2}w_{4r-2}+\cdots +w_{2r-2}w_{2r+2}$.
\end{enumerate}

\noindent Since the $Sq^{k}$ action on $H^{\ast }(B_{SO(n)};\mathbb{Z}_{2})$
is determined by the Wu--formula, the formula (2.2), together with the
algorithms given by Corollary 2.2, suffices to evaluate the $\gamma $ action
on $H^{\ast }(B_{SO(n)};\mathbb{Z}_{2})$. For example when $u=w_{2}$ we get
that

\begin{quote}
$w_{2}^{(1)}=w_{4},$

$w_{2}^{(2)}=w_{8}+w_{2}w_{6},$

$%
w_{2}^{(3)}=w_{16}+w_{2}w_{14}+w_{4}w_{12}+w_{6}w_{10}+w_{2}w_{6}w_{8}+w_{4}w_{6}^{2}+w_{2}w_{7}^{2} 
$

$\qquad
+w_{3}^{2}(w_{10}+w_{2}w_{8}+w_{4}w_{6})+w_{2}^{2}(w_{12}+w_{2}w_{10}+w_{4}w_{8}), 
$
\end{quote}

\noindent and in general, if $2^{k}\leq n$, that

\begin{quote}
$w_{2}^{(k)}=w_{2^{k}}+w_{2}w_{2^{k}-2}+\cdots +w_{2^{k-1}-2}w_{2^{k-1}+2}+$
higher terms.
\end{quote}

\noindent In contrast to the structure of the algebra $H^{\ast }(B_{SO(n)};%
\mathbb{Z}_{2})$ as a module over the Steenrod algebra \cite[Lemma 3.2]{PW},
these calculation reveal a striking property about the operator $\gamma $:
modulo the decomposable elements the derived sequence of $w_{2}$ consists of
all the $2$--power Stiefel--Whitney classes:

\begin{quote}
$\{w_{2}^{(0)},w_{2}^{(1)},w_{2}^{(2)},\cdots \}\equiv \{w_{2},w_{4},\cdots
,w_{2^{l(n)}},0,\cdots \}$, $l(n)=\left[ \ln n\right] $.
\end{quote}

\noindent This sequence will play a central role in the construction and
calculation of this paper.$\square $

\bigskip

An essential point of the operator $\alpha $ is that it always admits an
integral lift by i) of Theorem A. That is, it can be factored into $\rho
_{4}\circ f_{\alpha }$ for some $f_{\alpha }:H^{2r}(X;\mathbb{Z}%
_{2})\rightarrow H^{4r}(X)$. In the case $X=B_{SO(n)}$ of our interest a
canonical choice of such an integral lift $f_{\alpha }$ can be made
explicitly. Recall from Brown \cite{Br} and Feshbach \cite{Fe} that the
integral cohomology ring of $B_{SO(n)}$ is

\begin{enumerate}
\item[(2.3)] $H^{\ast }(B_{SO(n)})=\left\{ 
\begin{tabular}{l}
$\mathbb{Z}[p_{1},p_{2},\cdots ,p_{\left[ \frac{n-1}{2}\right]
},e_{n}]\oplus \tau (B_{SO(n)})$ if $n$ is even; \\ 
$\mathbb{Z}[p_{1},p_{2},\cdots ,p_{\left[ \frac{n-1}{2}\right] }]\oplus \tau
(B_{SO(n)})$ if $n$ is odd,%
\end{tabular}%
\right. $
\end{enumerate}

\noindent where $\tau (X)$ denotes the torsion ideal of the integral
cohomology $H^{\ast }(X)$ of a complex $X$, $2\tau (B_{SO(n)})=0$, and where 
$p_{i}$ (resp. $e_{n}$) is the $i^{th}$ Pontryagin class (resp. the Euler
class) of the canonical real $n$--bundle on $B_{SO(n)}$. With respect to
(2.3) Thomas introduced in \cite[\S 3]{Th} an operator $f:H^{r}(B_{SO(n)};%
\mathbb{Z}_{2})\rightarrow H^{2r}(B_{SO(n)})$, i.e. \textsl{the integral
representation}, by the following practical rules:

\begin{enumerate}
\item[(2.4)] $f(u):=\left\{ 
\begin{tabular}{l}
$p_{r}$, $\delta _{2}(Sq^{2r}w_{2r+1})$ or $e_{n}^{2}$ if $u=w_{2r},w_{2r+1}$
or $w_{n}$ if $n$ is even; \\ 
$f(w_{i_{1}})\cdots f(w_{i_{k}})$ if $u=w_{i_{1}}\cdots w_{i_{k}}$ is a
monomial; \\ 
$f(u_{1})+\cdots +f(u_{k})$ if $u=$ $u_{1}+\cdots +u_{k}$,%
\end{tabular}%
\right. $
\end{enumerate}

\noindent where $2r,2r+1<\left[ \frac{n-1}{2}\right] $, and where the $u_{i}$%
's are distinct monomials in $w_{2},\cdots ,w_{n}$. Based on Theorem A we
shall show that

\bigskip

\noindent \textbf{Theorem B.} \textsl{The pair }$(f,\gamma )$\textsl{\ of
operators on }$H^{\ast }(B_{SO(n)};\mathbb{Z}_{2})$ \textsl{satisfies,} 
\textsl{for any }$u\in S_{2}^{\ast }(B_{SO(n)})$\textsl{, that}

\begin{quote}
\textsl{i)} $\mathcal{B}(u)=\rho _{4}(f(u))+\theta (\gamma (u))$\textsl{,}

\textsl{ii)} $Sq^{1}(\gamma (u))=Sq^{2r}Sq^{1}(u)+u\cup Sq^{1}(u)$.
\end{quote}

\noindent \textbf{Example 2.5.} For $u=w_{2r}\in S_{2}^{\ast }(B_{SO(n)})$
we have $f(w_{2r})=p_{r}$ by the definition of $f$. Substituting this and
(2.2) into i) of Theorem B we obtain the formula

\begin{quote}
$\mathcal{B}(w_{2r})=\rho _{4}(p_{r})+\theta (w_{4r}+w_{2}w_{4r-2}+\cdots
+w_{2r-2}w_{2r+2})$
\end{quote}

\noindent implying that the $\func{mod}4$ reductions of the Pontryagin
classes of a manifold are homotopy invariants. This formula was originally
obtained by W. Wu \cite{Wu} by computing with the Schubert cells on $%
B_{SO(n)}$. S.S. Chern suggested a different approach which was implemented
by Thomas in \cite[Theorem C]{Th1}.

Concerning this topic property ii) of Theorem A may be called \textsl{the} 
\textsl{generalized Wu--formula }on $\delta _{2}$--formal spaces.$\square $

\bigskip

Turning to our main concerns the classifying spaces $B_{Spin(n)}$ and $%
B_{Spin^{c}(n)}$ fit into the fibered sequences

\begin{enumerate}
\item[(2.5)] $\mathbb{C}P^{\infty }\overset{i}{\rightarrow }B_{Spin^{c}(n)}%
\overset{\pi }{\rightarrow }B_{SO(n)}$ and

\item[(2.6)] $U(1)\rightarrow B_{Spin(n)}\overset{\psi }{\rightarrow }%
B_{Spin^{c}(n)}\overset{\iota }{\rightarrow }\mathbb{C}P^{\infty }$,
\end{enumerate}

\noindent where the maps $\pi $ and $\iota $ are induced by the obvious
epimorphisms

\begin{quote}
$Spin(n)\times _{\mathbb{Z}_{2}}U(1)\rightarrow SO(n)$ and $Spin(n)\times _{%
\mathbb{Z}_{2}}U(1)\rightarrow U(1)$,
\end{quote}

\noindent respectively. Let $\{w_{2},w_{2}^{(1)},w_{2}^{(2)},\cdots \}$ be
the derived sequence of the second Stiefel Whitney class $w_{2}$. Applying
the operator $f$ gives rise to the sequence $\{f(w_{2}),$ $%
f(w_{2}^{(1)}),f(w_{2}^{(2)}),\cdots \}$ of integral cohomology classes on $%
B_{SO(n)}$. By examining the $\pi ^{\ast }$ images of these two sequences in
the cohomology of $B_{Spin^{c}(n)}$ we single out a set of new generators of
the ring $H^{\ast }(B_{Spin^{c}(n)})$ in the following result, where $x$\
denotes the Euler class of the Hopf line bundle $\lambda $\ on $\mathbb{C}%
P^{\infty }$.

\bigskip

\noindent \textbf{Theorem C. }\textsl{There exists a unique sequence }$%
\left\{ q_{r}\QTR{sl}{,\ }r\geq 0\right\} $ \textsl{of} \textsl{integral} 
\textsl{cohomology classes on }$B_{Spin^{c}(n)}$\textsl{, }$\deg
q_{r}=2^{r+1}$\textsl{,} \textsl{that satisfies the following system}

\begin{quote}
\textsl{i)} $q_{0}=\iota ^{\ast }(x)$\textsl{,}

\textsl{ii)} $\rho _{2}(q_{r})=\pi ^{\ast }(w_{2}^{(r)})$\textsl{; }

\textsl{iii)} $2q_{r+1}+q_{r}^{2}=\pi ^{\ast }f(w_{2}^{(r)})$\textsl{, }$%
r\geq 0$\textsl{.}
\end{quote}

For an integer $n\geq 7$ (see Remark 1.1) we set $h(n)=\left[ \frac{n-1}{2}%
\right] $, and let $\theta _{n}\in H^{\ast }(B_{Spin^{c}(n)})$ be the Euler
class of the complex bundle $\xi _{n}$ associated to the complex spin
presentation $Spin^{c}(n)\rightarrow U(2^{h(n)})$ (\cite{ABS,HK}). Regarding
the cohomology $H^{\ast }(B_{Spin^{c}(n)})$ as a module over its subring $%
\pi ^{\ast }H^{\ast }(B_{SO(n)})$, our main result presents the cohomology $%
H^{\ast }(B_{Spin^{c}(n)})$ by the unique sequence $\left\{ q_{r},r\geq
0\right\} $ obtained by Theorem C, together with the Euler class $\theta
_{n} $.

\bigskip

\noindent \textbf{Theorem D. }\textsl{The cohomology of }$B_{Spin^{c}(n)}$%
\textsl{\ has the presentation}

\begin{enumerate}
\item[(2.6)] $H^{\ast }(B_{Spin^{c}(n)})=\pi ^{\ast }H^{\ast
}(B_{SO(n)})\otimes \mathbb{Z}[q_{0},q_{1},\cdots ,q_{h(n)-1},\theta
_{n}]/R_{n}$,
\end{enumerate}

\noindent \textsl{in which }$R_{n}$ \textsl{denotes the ideal generated by
the following elements}

\begin{quote}
\textsl{i)} $2q_{r+1}+q_{r}^{2}-\pi ^{\ast }f(w_{2}^{(r)})$\textsl{,\quad }$%
0\leq r\leq h(n)-2\QTR{sl}{;}$

\textsl{ii)} $\pi ^{\ast }\delta _{2}(z)\cup q_{r}-\pi ^{\ast }\delta
_{2}(z\cup w_{2}^{(r)}),$ $0\leq r\leq h(n)-1$\textsl{;}

\textsl{iii)} $4(-1)^{h(n)}\theta _{n}+q_{h(n)-1}^{2}-a_{n}$\textsl{,}
\end{quote}

\noindent \textsl{where} $z\in H^{\ast }(B_{SO(n)},\mathbb{Z}_{2})$\textsl{,
and where }$a_{n}\in \pi ^{\ast }H^{+}(B_{SO(n)})\otimes \mathbb{Z}%
(q_{0},q_{1},\cdots ,q_{h(n)-1})$ \textsl{is an} \textsl{element to be
specified in Lemma 5.3.}

\bigskip

Leaving out $\otimes $ for notational simplicity the relations i) and ii) of
Theorem D are inherited from the properties iii) and ii) of Theorem C. They
express the relationship between the two sequences $\left\{ q_{r},r\geq
0\right\} $ and $\left\{ \pi ^{\ast }p_{r},r\geq 1\right\} $ of integral
cohomology classes on $B_{Spin^{c}(n)}$, and deal with the product between
the free part and the torsion ideal $\tau (B_{Spin^{c}(n)})$ of the ring $%
H^{\ast }(B_{Spin^{c}(n)})$, respectively.

An outline of the paper is as follows. Sections \S 3 to \S 5 are devoted to
show Theorems A--D. The calculation is extended in Section \S 6 to obtain a
similar presentation of the ring $H^{\ast }(B_{Spin(n)})$ in Theorem D'. The
remaining sections constitute applications of Theorem D'. We determine in \S %
8 the ring of integral Weyl invariants of the group $Spin(n)$; introduce in 
\S 9 the spin characteristic classes for the spin vector bundles. In spin
geometry the Spin characteristic classes can play roles that may not be
replaced by the regular characteristic classes. Section \S 10 is devoted to
present such examples.

Effective computability with characteristic classes is essential in geometry
and invariant theory. To illustrate the algorithmic nature of the operators
and results developed in this paper a number of examples are included. In
particular, in a concrete situation Theorems D and D' are directly
applicable to present the rings $H^{\ast }(B_{Spin^{c}(n)})$ and $H^{\ast
}(B_{Spin(n)})$ by a minimal system of explicit generators and relations,
see in Examples 5.4 and 6.5; Theorems 7.7 and 7.8 are ready to deduce the
integral Weyl invariants of the groups $Spin^{c}(n)$ and $Spin(n)$, which
are shown in Examples 7.6 and 7.9.

\section{The cohomology operation $F=(\protect\gamma ,\protect\alpha )$}

For a $\func{mod}2$ cohomology class $u\in H^{2r}(X;\mathbb{Z}_{2})$ of a CW
complex $X$ the Pontryagin square $\mathcal{B}(u)\in H^{4r}(X;\mathbb{Z}%
_{4}) $ can be defined by the formula

\begin{quote}
$\mathcal{B}(u)\equiv \rho _{4}(\widetilde{u}\cup _{0}\widetilde{u}+\delta (%
\widetilde{u})\cup _{1}\widetilde{u})$ (see \cite{BT}),
\end{quote}

\noindent where $\widetilde{u}$ is an integral lift of $u$ in the cochain
complex $C^{\ast }(X;\mathbb{Z})$ associated to $X$, and $\cup _{i}$ denotes
the $i^{th}$ cup product on $C^{\ast }(X;\mathbb{Z})$. Based on this formula
a co--chain level calculation verifies the following universal relations

\begin{enumerate}
\item[(3.1)] $\delta _{2}(u\cup u)=2\delta _{4}(\mathcal{B}(u))$ in $%
H^{4r+1}(X)$;

\item[(3.2)] $\rho _{2}\delta _{4}\mathcal{B}(u)=Sq^{2r}Sq^{1}u+u\cup
Sq^{1}u $ in $H^{4r+1}(X;\mathbb{Z}_{2})$.
\end{enumerate}

\noindent \textbf{Proof of Theorem A.} Let $X$ be a $\delta _{2}$ formal
space. For each $u\in H^{2r}(X;\mathbb{Z}_{2})$ the property $\delta
_{2}(u\cup u)=0$ implies that $\delta _{4}(\mathcal{B}(u))\in \func{Im}%
\delta _{2}$ by (3.1). In view of the isomorphism $\delta _{2}:S_{2}^{\ast
}(X)\cong \func{Im}\delta _{2}$ there exists a unique $u_{1}\in
S_{2}^{4r}(X) $ so that

\begin{enumerate}
\item[(3.3)] $\delta _{2}(u_{1})=\delta _{4}(\mathcal{B}(u))$.
\end{enumerate}

\noindent We can now formulate the desired operation $F=(\alpha ,\gamma
):H^{2r}(X;\mathbb{Z}_{2})\rightarrow H^{4r}(X;\mathbb{Z}_{4})\times
S^{4r}(X;\mathbb{Z}_{2})$ by setting

\begin{quote}
$\gamma (u):=u_{1}$ and $\alpha (u):=\mathcal{B}(u)-\theta (\gamma (u))$.
\end{quote}

Applying $\rho _{2}$ to both sides of (3.3) we get by (3.2) that

\begin{quote}
$Sq^{1}\gamma (u)=Sq^{2r}Sq^{1}u+u\cup Sq^{1}u$,
\end{quote}

\noindent showing property iii) of Theorem A. From $\delta _{4}\circ \theta
=\delta _{2}$ and (3.3) we obtain

\begin{quote}
$\delta _{4}\alpha (u)=\delta _{4}(\mathcal{B}(u)-\theta (\gamma
(u)))=\delta _{4}(\mathcal{B}(u))-\delta _{2}(\gamma (u))=0$,
\end{quote}

\noindent implying $\alpha (u)\in \func{Im}\rho _{4}$ (i.e. property i) of
Theorem A).

Summarizing, the pair $F=(\gamma ,\alpha )$ of operators fulfills the
properties i), ii) and iii) of Theorem A, whose uniqueness comes directly
from its definition.$\square $

\bigskip

\noindent \textbf{Proof of Theorem B.} Let $f:H^{\ast }(B_{SO(n)};\mathbb{Z}%
_{2})\rightarrow H^{\ast }(B_{SO(n)})$ be the map entailed in (2.4). It has
been shown by Thomas \cite[Lemma (3.9)]{Th} that for any $u\in S_{2}^{\ast
}(B_{SO(n)})$ there exists a unique element $v\in S_{2}^{\ast }(B_{SO(n)})$
so that

\begin{enumerate}
\item[(3.5)] $\mathcal{B}(u)=\rho _{4}(f(u))+\theta (v)$.
\end{enumerate}

\noindent It suffices for us to show that $v=\gamma (u)$. Since the space $%
B_{SO(n)}$ is $\delta _{2}$--formal applying $\delta _{4}$ to both sides of
(3.5) we get by (3.3) that

\begin{quote}
$\delta _{2}(\gamma (u))=\delta _{4}(\theta (v))=\delta _{2}(v)$ (since $%
\delta _{4}\circ \theta =\delta _{2}$).
\end{quote}

\noindent With $\gamma (u),v\in S_{2}^{\ast }(B_{SO(n)})$ while $\delta _{2}$
injects on $S_{2}^{\ast }(B_{SO(n)})$, we obtain $v=$ $\gamma (u)$.$\square $

\section{The proof of Theorem C}

Applying the Serre spectral sequence to the fibration

\begin{enumerate}
\item[(4.1)] $\mathbb{C}P^{\infty }\overset{i}{\rightarrow }B_{Spin^{c}(n)}%
\overset{\pi }{\rightarrow }B_{SO(n)}$ (see (2.5))
\end{enumerate}

\noindent the cohomologies of the total space $B_{Spin^{c}(n)}$ with fields
coefficients have been computed. We begin by recalling the relevant results
due to Borel, Hirzebruch, Quillen, Harada and Kono. Let $\mathbb{Z}_{0}$ be
the field of rationals and denote by $\rho _{0}$ be the cohomology
homomorphism induced by the inclusion $\mathbb{Z}\subset \mathbb{Z}_{0}$.
Set $q_{0}:=\iota ^{\ast }(x)$, where $x$ is the Euler class of the Hopf
line bundle $\lambda $ on $\mathbb{C}P^{\infty }$.

\bigskip

\noindent \textbf{Lemma 4.1.} \textsl{If either }$p=0$\textsl{\ or }$p\geq 3$%
\textsl{\ is a prime, then }$H^{\ast }(B_{Spin^{c}(n)};\mathbb{Z}_{p})$%
\textsl{\ is a free polynomial algebra on the generators}

\begin{enumerate}
\item[(4.2)] $\rho _{p}(q_{0}),\rho _{p}(\pi ^{\ast }p_{1}),\cdots ,\rho
_{p}(\pi ^{\ast }p_{\left[ \frac{n-1}{2}\right] }),$ \textsl{and} $\rho
_{p}(\pi ^{\ast }e_{n})$ \textsl{if} $n\equiv 0\func{mod}2$.
\end{enumerate}

\textsl{In addition, the map }$\rho :H^{\ast }(B_{Spin^{c}(n)})\rightarrow
H^{\ast }(B_{Spin^{c}(n)};\mathbb{Z}_{0})\times H^{\ast }(B_{Spin^{c}(n)};%
\mathbb{Z}_{2})$ \textsl{by} $\rho (z)=(\rho _{0}(z),\rho _{2}(z))$\textsl{\
injects.}

\bigskip

\noindent \textbf{Proof.} Since the composition $\iota \circ i:\mathbb{C}%
P^{\infty }\rightarrow \mathbb{C}P^{\infty }$ (see (2.5) and (2.6)) is of
degree $2$ the class $i^{\ast }(q_{0})=2x$ generates the algebra $H^{\ast }(%
\mathbb{C}P^{\infty };\mathbb{Z}_{p})=\mathbb{Z}_{p}[x]$ for all $p\neq 2$.
The first assertion follows from the Leray--Hirsch property \cite[p.231]{Hus}
of the fibration (4.1) with $\mathbb{Z}_{p}$ coefficients.

According to Borel and Hirzebruch \cite[30.6.]{BH} if $X$ is a space with $%
2\tau (X)=0$, then the map $\rho :H^{\ast }(X)\rightarrow H^{\ast }(X;%
\mathbb{Z}_{0})\times H^{\ast }(X;\mathbb{Z}_{2})$ by $\rho (z)=(\rho
_{0}(z),\rho _{2}(z))$ injects. The second assertion follows from $2\tau
(B_{Spin^{c}(n)})=0$ by Harada and Kono \cite[Theorem 3.7]{HK}.$\square $

\bigskip

Turning to the algebra $H^{\ast }(B_{Spin^{c}(n)};\mathbb{Z}_{2})$ the
transgression $\sigma $ in the fibration (4.1) clearly satisfies $\sigma
(\rho _{2}(x))=w_{3}$. Since $\sigma $ commutes with the Steenrod squares
the standard relation $Sq^{2^{k}}x_{k}=x_{k+1}$ on $H^{\ast }(\mathbb{C}%
P^{\infty };\mathbb{Z}_{2})$ implies that

\begin{enumerate}
\item[(4.3)] $\sigma (x_{k+1})=Sq^{2^{k}}\sigma (x_{k})$, $k\geq 0$, where $%
x_{k}:=(\rho _{2}(x))^{2^{k-1}}$.
\end{enumerate}

\noindent For a subset $\{a_{1},\cdots ,a_{r}\}$ of an algebra $A$ denote by 
$\left\langle a_{1},\cdots ,a_{r}\right\rangle $ the ideal generated by $%
a_{1},\cdots ,a_{r}$, and let $A/\left\langle a_{1},\cdots
,a_{r}\right\rangle $ be the quotient algebra.

\bigskip

\noindent \textbf{Lemma 4.2.} \textsl{Let }$J=\left\langle \sigma
(x_{1}),\cdots ,\sigma (x_{h(n)})\right\rangle $\textsl{, }$h(n)=\left[ 
\frac{n-1}{2}\right] $\textsl{, and let }$\theta _{n}$\textsl{\ be the Euler
class of the presentation }$Spin^{c}(n)\rightarrow U(2^{h(n)})$\textsl{. Then%
}

\begin{enumerate}
\item[(4.4)] $H^{\ast }(B_{Spin^{c}(n)};\mathbb{Z}_{2})=H^{\ast }(B_{SO(n)};%
\mathbb{Z}_{2})/J\otimes \mathbb{Z}_{2}[\theta _{n}]$\textsl{.}
\end{enumerate}

\textsl{In particular, the torsion ideal }$\tau (B_{Spin^{c}(n)})$\textsl{\
of the ring }$H^{\ast }(B_{Spin^{c}(n)})$\textsl{\ is}

\begin{enumerate}
\item[(4.5)] $\tau (B_{Spin^{c}(n)})=\pi ^{\ast }\tau (B_{SO(n)})\otimes 
\mathbb{Z}[\theta _{n}]$.
\end{enumerate}

\noindent \textbf{Proof.} Formula (4.4) goes to Harada and Kono \cite[%
Theorem 3.5]{HK} (see also Quillen \cite[Theorem 6.5]{Q}). With $2\tau
(B_{SO(n)})=2\tau (B_{Spin^{c}(n)})=0$ the map $\pi $ induces the
commutative diagram

\begin{quote}
\begin{tabular}{llll}
$0\rightarrow $ & $\tau (B_{SO(n)})$ & $\overset{\rho _{2}}{\rightarrow }$ & 
$H^{\ast }(B_{SO(n)};\mathbb{Z}_{2})$ \\ 
& $\pi ^{\ast }\downarrow $ &  & $\pi ^{\ast }\downarrow $ \\ 
$0\rightarrow $ & $\tau (B_{Spin^{c}(n)})$ & $\overset{\rho _{2}}{%
\rightarrow }$ & $H^{\ast }(B_{Spin^{c}(n)};\mathbb{Z}_{2})$%
\end{tabular}
\end{quote}

\noindent in which both $\rho _{2}$ inject. Since $\tau (B_{Spin^{c}(n)})$
is an ideal the map

\begin{quote}
$h:\pi ^{\ast }\tau (B_{SO(n)})\otimes \mathbb{Z}[\theta _{n}]\rightarrow
\tau (B_{Spin^{c}(n)})$ by $h(\pi ^{\ast }x\otimes \theta _{n})=\pi ^{\ast
}x\cup \theta _{n}$
\end{quote}

\noindent is well defined, and gives rise to the isomorphism (4.5).

Indeed, by (4.4) the composition $\rho _{2}\circ h$ injects, hence $h$
injects, too. On the other hand for any $x\in \tau (B_{Spin^{c}(n)})$ there
exists $y\in H^{\ast }(B_{Spin^{c}(n)};\mathbb{Z}_{2})$ so that $\delta
_{2}(y)=x$. By formula (4.4) the map $h$ also surjects.$\square $

\bigskip

As the space $B_{SO(n)}$ is $\delta _{2}$--formal the derived sequence $%
\{w_{2},w_{2}^{(1)},\cdots \}$ of the Stiefel Whitney class $w_{2}$ is
defined, see Example 2.4. Its relationship with the sequence $\{\sigma
(x_{1}),\sigma (x_{2}),\cdots \}$ defined by (4.3) is stated in the
following result.

\bigskip

\noindent \textbf{Lemma 4.3.} \textsl{In }$H^{\ast }(B_{SO(n)};\mathbb{Z}%
_{2})$\textsl{\ let }$J_{n,k}=$\textsl{\ }$\left\langle \sigma
(x_{1}),\cdots ,\sigma (x_{k})\right\rangle $\textsl{, }$k\geq 0$\textsl{.
Then}

\begin{enumerate}
\item[(4.6)] $Sq^{1}(w_{2}^{(k-1)})=\sigma (x_{k})+\beta _{k-1}$ \textsl{%
with }$\beta _{k-1}\in J_{n,k-1}$\textsl{, }$k\geq 1$\textsl{.}
\end{enumerate}

\noindent \textbf{Proof.} If $k=1$ the formula (4.6) is verified by $\sigma
(x_{1})=w_{3}=Sq^{1}(w_{2})$. Assume next that it holds for some $k\geq 1$.
Then

\begin{quote}
$Sq^{1}(w_{2}^{(k)})=Sq^{1}(\gamma (w_{2}^{(k-1)}))$ (by $w_{2}^{(k)}=\gamma
(w_{2}^{(k-1)})$)

$=Sq^{2^{k}}Sq^{1}(w_{2}^{(k-1)})+w_{2}^{(k-1)}\cup Sq^{1}(w_{2}^{(k-1)})$
(by iii) of Theorem A)

$=Sq^{2^{k}}(\sigma (x_{k})+\beta _{k-1})+w_{2}^{(k-1)}\cup (\sigma
(x_{k})+\beta _{k-1})$ (by induction)

$=\sigma (x_{k+1})+Sq^{2^{k}}\beta _{k-1}+w_{2}^{(k-1)}\cup (\sigma
(x_{k})+\beta _{k-1})$ (by (4.3)).
\end{quote}

\noindent The inductive procedure showing (4.6) is completed by taking

\begin{quote}
$\beta _{k}:=Sq^{2^{k}}\beta _{k-1}+w_{2}^{(k-1)}\cup (\sigma (x_{k})+\beta
_{k-1})$.$\square $
\end{quote}

We proceed to a constructive proof of Theorem C. Since $\pi ^{\ast }\circ
\sigma =0$ we get by (4.5) that $Sq^{1}(\pi ^{\ast }(w_{2}^{(k)}))=0$. With
the space $B_{Spin^{c}(n)}$ being $\delta _{2}$--formal we obtain further $%
\delta _{2}(\pi ^{\ast }(w_{2}^{(k)}))=0$, $k\geq 0$. In view of the
Bockstein exact sequence

\begin{center}
$\cdots \rightarrow H^{r}(B_{Spin^{c}(n)})\overset{\rho _{2}}{\rightarrow }%
H^{r}(B_{Spin^{c}(n)};\mathbb{Z}_{2})\overset{\delta _{2}}{\rightarrow }%
H^{r+1}(B_{Spin^{c}(n)})\rightarrow \cdots $
\end{center}

\noindent this implies that the classes $\pi ^{\ast }(w_{2}^{(k)})$ admit
integral lifts

\begin{enumerate}
\item[(4.7)] $\rho _{2}(q_{k}^{\prime })=\pi ^{\ast }(w_{2}^{(k)})$ for some 
$q_{k}^{\prime }\in H^{2^{k+1}}(B_{Spin^{c}(n)})$, $k\geq 0$.
\end{enumerate}

\noindent In particular,

\begin{quote}
$\mathcal{B}(\pi ^{\ast }(w_{2}^{(k)}))=\rho _{4}(q_{k}^{\prime }\cup
q_{k}^{\prime })$, $\theta (\pi ^{\ast }(w_{2}^{(k)}))=$ $2\rho
_{4}(q_{k}^{\prime })$.
\end{quote}

\noindent Thus, applying $\pi ^{\ast }$ to $u=w_{2}^{(k)}$ the relation i)
of Theorem B becomes

\begin{quote}
$\rho _{4}(q_{k}^{\prime }\cup q_{k}^{\prime })=\rho _{4}(\pi ^{\ast
}f(w_{2}^{(k)}))+2\rho _{4}(q_{k+1}^{\prime })$ on $H^{\ast
}(B_{Spin^{c}(n)};\mathbb{Z}_{4})$,
\end{quote}

\noindent implying that there exist integral class $v_{k+1}\in H^{\ast
}(B_{Spin^{c}(n)})$ so that

\begin{enumerate}
\item[(4.8)] $2q_{k+1}^{\prime }+q_{k}^{\prime }\cup q_{k}^{\prime }=\pi
^{\ast }f(w_{2}^{(k)})+4v_{k+1}$.
\end{enumerate}

\noindent \textbf{Proof of Theorem C. }Since the class $q_{0}=\iota ^{\ast
}(x)$ generates $H^{2}(B_{Spin^{c}(n)})=\mathbb{Z}$ with $\rho
_{2}(q_{0})=\pi ^{\ast }(w_{2})$ we can take in (4.7) that $q_{0}^{\prime
}=q_{0}$, and define in term of (4.8) that $q_{1}:=q_{1}^{\prime }-2v_{1}$.
Then, the relations ii) and iii) of Theorem C for the case $r=1$ are
verified by

\begin{quote}
$\rho _{2}(q_{1})=\rho _{2}(q_{1}^{\prime })=\pi ^{\ast }(w_{2}^{(1)})$ (by
(4.7));

$2q_{1}+q_{0}\cup q_{0}=\pi ^{\ast }f(w_{2}^{(0)})$ (by (4.8)).
\end{quote}

Assume next that a sequence $q_{0},\cdots ,q_{r}$ of classes satisfying the
properties i), ii) and iii) of Theorem C has been obtained for some $r\geq 1$%
. Take in (4.7) that $q_{r}^{\prime }=q_{r}$ and define in term of (4.8)
that $q_{r+1}:=q_{r+1}^{\prime }-2v_{r+1}$. Then

\begin{quote}
$\rho _{2}(q_{r+1})=\pi ^{\ast }(w_{2}^{(r+1)})$, $2q_{r+1}+q_{r}\cup
q_{r}=\pi ^{\ast }f(w_{2}^{(r)})$.
\end{quote}

\noindent This completes the inductive construction of a sequence $%
\{q_{r},r\geq 0\}$ fulfilling the system i)--iii) of Theorem C.

To see the uniqueness of the sequence $\{q_{r},r\geq 0\}$ we make use of the
injection $\rho $ in Lemma 4.1. Note that the properties i), ii) and iii) of
Theorem C suffices to decide the $\rho $--image of $q_{r}$ as $\rho
(q_{r})=(g_{r},\pi ^{\ast }w_{2}^{(r)})$, $r\geq 1$, where $g_{r}\in H^{\ast
}(B_{Spin^{c}(n)};\mathbb{Z}_{0})$ is the unique polynomial (with rational
coefficients) in the generators (4.2) defined recurrently by the relation
iii) of Theorem C as

\begin{quote}
$\rho _{0}(q_{1})=g_{1}:=\frac{1}{2}(\rho _{0}(\pi ^{\ast }p_{1})-\rho
_{0}(q_{0})^{2})$

$\rho _{0}(q_{2})=g_{2}:=\frac{1}{2}(\rho _{0}(\pi ^{\ast }p_{2})-g_{1}^{2})$%
, $\cdots $,

$\rho _{0}(q_{r})=g_{r}:=\frac{1}{2}(\rho _{0}(\pi ^{\ast
}f(w_{2}^{(r-1)}))-g_{r-1}^{2}))$, $r\geq 2$.
\end{quote}

\noindent The injectivity of $\rho $ implies that, if $\left\{ q_{r}^{\prime
},1\leq r\right\} $ is a second sequence satisfying the system i)--iii) of
Theorem C, then $q_{r}^{\prime }=q_{r}$, as required.$\square $

\bigskip

We conclude this section with two applications of Theorem C. Firstly, let $%
H^{+}(B_{SO(n)})$ be the subring of $H^{\ast }(B_{SO(n)})$ consisting of
elements in the positive degrees. The induced action of the fiber inclusion $%
i$ on cohomology is determined in the following result.

\bigskip

\noindent \textbf{Lemma 4.4. }\textsl{The map }$i^{\ast }$\textsl{\
satisfies that }$i^{\ast }\circ \pi ^{\ast }=0$\textsl{\ on }$%
H^{+}(B_{SO(n)})$\textsl{, and that}

\begin{enumerate}
\item[(4.9)] $i^{\ast }(q_{r})=2(-1)^{r}x^{2^{r}}$, $r\geq 0$; $\quad
i^{\ast }(\theta _{n})=x^{2^{h(n)}}$.
\end{enumerate}

\noindent \textbf{Proof.} As the composition $\pi \circ i$ is null--homotopy
we get $i^{\ast }\circ \pi ^{\ast }=0$. Consequently, applying $i^{\ast }$
to the relation iii) of Theorem C yields

\begin{quote}
$2i^{\ast }(q_{r+1})+i^{\ast }(q_{r})^{2}=0$, $r\geq 0$.
\end{quote}

\noindent Inputting $i^{\ast }(q_{0})=2x$ we obtain the first relation in
(4.9) by induction on $r$. The second one is verified by the geometric fact $%
i^{\ast }\xi _{n}=\lambda \oplus \cdots \oplus \lambda $ ($2^{h(n)}$ copies).%
$\square $

\bigskip

Recall that the $\func{mod}2$ Bockstein cohomology $H_{\beta }^{\ast }(X)$
of a space $X$ is the kernel modulo the image of the operation $\beta
=Sq^{1} $ on $H^{\ast }(X;\mathbb{Z}_{2})$. Granted with Theorem C we can
express the Bockstein $H_{\beta }^{\ast }(B_{Spin^{c}(n)})$ by the mod $2$
reductions of explicit integral cohomology classes on $B_{Spin^{c}(n)}$. For
a set $\{y_{1},\cdots ,y_{r}\}$ of graded elements let $\Delta (y_{1},\cdots
,y_{r})$ be the graded free $\mathbb{Z}$ module with the basis

\begin{quote}
$\{1,y_{i_{1}}y_{i_{2}}\cdots y_{i_{k}}$, $1\leq i_{1}<\cdots <i_{k}\leq r\}$
\end{quote}

\noindent (i.e. the square free monomials in the $y_{1},\cdots ,y_{r}$). For
a prime $p$ we use the notion $\Delta _{p}(y_{1},\cdots ,y_{r})$ to simplify
the tensoring $\mathbb{Z}_{p}\otimes \Delta (y_{1},\cdots ,y_{r})$.

\bigskip

\noindent \textbf{Lemma 4.5.} \textsl{Let }$\{q_{r},r\geq 0\}$ \textsl{be
the unique sequence of integral cohomology classes on }$B_{Spin^{c}(n)}$\ 
\textsl{obtained in Theorem C. Then}

\begin{enumerate}
\item[(4.10)] $H_{\beta }^{\ast }(B_{Spin^{c}(n)})=\pi ^{\ast }H_{\beta
}^{\ast }(B_{SO(n)})\otimes \Delta _{2}(q_{0},\cdots ,q_{h(n)-1})\otimes 
\mathbb{Z}_{2}[\theta _{n}]$\textsl{,}
\end{enumerate}

\noindent \textsl{where, for }$n=2k+1$\textsl{\ or }$2k$\textsl{,}

\begin{quote}
$\pi ^{\ast }H_{\beta }^{\ast }(B_{SO(n)})=\mathbb{Z}_{2}[w_{2}^{2},\cdots
,w_{2k}^{2}]$ \textsl{or} $\mathbb{Z}_{2}[w_{2}^{2},\cdots
,w_{2(k-1)}^{2},w_{2k}]$\textsl{,}
\end{quote}

\noindent \textsl{and where, with }$\mathbb{Z}_{2}$ \textsl{coefficients we
can write }$q_{r}$\textsl{, }$\theta _{n}$\textsl{\ in places of }$\rho
_{2}(q_{r})$\textsl{, }$\rho _{2}(\theta _{n})$\textsl{.}

\bigskip

\noindent \textbf{Proof.} A sequence $\{a_{1},\cdots ,a_{r}\}$ of elements
in the algebra $H^{\ast }(B_{SO(n)};\mathbb{Z}_{2})$ is called \textsl{%
regular} if for each $i\geq 1$ the residue class $\left[ a_{i}\right] $ in
the quotient algebra

\begin{quote}
$H^{\ast }(B_{SO(n)};\mathbb{Z}_{2})/\left\langle a_{1},\cdots
,a_{i-1}\right\rangle $
\end{quote}

\noindent is not a zero divisor. According to Quillen \cite{Q} the sequence $%
\left\{ \sigma (x_{1}),\cdots ,\sigma (x_{h(n)-1})\right\} $ is regular in $%
H^{\ast }(B_{SO(n)};\mathbb{Z}_{2})$ (see Remark 4.6 below). From formula
(4.6) we can conclude further that

a) the sequence $\{Sq^{1}(w_{2}^{(0)}),\cdots ,Sq^{1}(w_{2}^{(h(n)-1)})\}$
is regular;

b) $J_{n,k}=I_{n,k}$, $1\leq k\leq h(n)-1$, where $I_{n,k}=\left\langle
Sq^{1}(w_{2}^{(0)}),\cdots ,Sq^{1}(w_{2}^{(k-1)})\right\rangle $.

\noindent Combining formula (4.4) with property b) we get the more useful
presentation

\begin{quote}
$H^{\ast }(B_{Spin^{c}(n)};\mathbb{Z}_{2})=H^{\ast }(B_{SO(n)};\mathbb{Z}%
_{2})/I_{n,h(n)}\otimes \mathbb{Z}_{2}[\theta _{n}]$,
\end{quote}

\noindent of the algebra $H^{\ast }(B_{Spin^{c}(n)};\mathbb{Z}_{2})$, on
which the operator $\beta $ is determined by

\begin{quote}
$\beta (\pi ^{\ast }(w_{2r}))=\pi ^{\ast }(w_{2r+1})$, $\beta (\pi ^{\ast
}(w_{2r+1}))=0$, $\beta (\rho _{2}(\theta _{n}))=0$,

$\beta (Sq^{1}(w_{2}^{(k)}))=Sq^{1}Sq^{1}w_{2}^{(k)}=0$, $0\leq k\leq h(n)-1$%
.
\end{quote}

\noindent As a result, if we introduce the quotient algebras $C_{k}=H^{\ast
}(B_{SO(n)};\mathbb{Z}_{2})/I_{n,k}$, $0\leq k\leq h(n)-1$, then the pairs $%
\left\{ C_{k},\beta \right\} $ become $\beta $--complexes that fit into the
short exact sequences

\begin{quote}
$0\rightarrow \left\{ C_{k},\beta \right\} \overset{\cup Sq^{1}(w_{2}^{(k)})}%
{\rightarrow }\left\{ C_{k},\beta \right\} \rightarrow \left\{ C_{k+1},\beta
\right\} \rightarrow 0,$
\end{quote}

\noindent where the first map $\cup Sq^{1}(w_{2}^{(i)})$ injects (because $%
Sq^{1}(w_{2}^{(k)})\in C_{k}$ is not a zero divisor by property a)), and
where the second map is the obvious quotient. Direct calculation shows that
its induced long exact sequences in cohomologies break into the short exact
sequences

\begin{quote}
$0\rightarrow H_{\beta }^{\ast }(C_{k})\rightarrow H_{\beta }^{\ast
}(C_{k+1})\rightarrow \rho _{2}(q_{k})\cup H_{\beta }^{\ast
}(C_{k})\rightarrow 0$,
\end{quote}

\noindent implying that $H_{\beta }^{\ast }(C_{k+1})$ is a free module over $%
H_{\beta }^{\ast }(C_{k})$ with basis $\{1,\rho _{2}(q_{k})\}$. Arguing
inductively we obtain the additive presentation

\begin{quote}
$H_{\beta }^{\ast }(C_{h(n)})=H_{\beta }^{\ast }(C_{0})\otimes \Delta
_{2}(q_{0},\cdots ,q_{h(n)-1})$
\end{quote}

\noindent in which $H_{\beta }^{\ast }(C_{0})=\pi ^{\ast }H_{\beta }^{\ast
}(B_{SO(n)})$ by the definition. Finally, we obtain (4.10) by the K\"{u}%
nneth formula

\begin{quote}
$H_{\beta }^{\ast }(C_{h(n)}\otimes \mathbb{Z}_{2}[\theta _{n}])=H_{\beta
}^{\ast }(C_{h(n)})\otimes \mathbb{Z}_{2}[\theta _{n}]$.$\square $
\end{quote}

\bigskip

\noindent \textbf{Remark 4.6. }In \cite[(6.2)]{Q} Quillen defined the
function $h$ on the positive integers whose values are given by the
following table

\begin{center}
\begin{tabular}{l||llllllll}
\hline
$n$ & $8l+1$ & $8l+2$ & $8l+3$ & $8l+4$ & $8l+5$ & $8l+6$ & $8l+7$ & $8l+8$
\\ \hline
$h(n)$ & $4l$ & $4l+1$ & $4l+2$ & $4l+2$ & $4l+3$ & $4l+3$ & $4l+3$ & $4l+3$
\\ \hline
\end{tabular}%
.
\end{center}

\noindent He stated in \cite[Theorem 6.3]{Q} that the sequence $w_{2},\sigma
(x_{1}),\cdots ,\sigma (x_{h(n)})$ is regular on $H^{\ast }(B_{SO(n)};%
\mathbb{Z}_{2})$. This statement should be corrected as

\begin{quote}
"the sequence $w_{2},\sigma (x_{1}),\cdots ,\sigma (x_{h(n)-1})$ is regular
in $H^{\ast }(B_{SO(n)};\mathbb{Z}_{2})$."
\end{quote}

\noindent Indeed, direct computation shows that

\begin{quote}
$\sigma (x_{1})=w_{3},$

$\sigma (x_{2})=Sq^{2}(\sigma (x_{1}))=w_{5}+w_{2}w_{3}$

$\sigma (x_{3})=Sq^{4}(\sigma
(x_{2}))=w_{9}+w_{2}w_{7}+w_{3}w_{6}+w_{3}^{3}+w_{2}^{2}(w_{5}+w_{2}w_{3})$
\end{quote}

\noindent confirming that \cite[Theorem 6.3]{Q} already fails for $n=5,6,7,8$%
.$\square $

\section{The cohomology $H^{\ast }(B_{Spin^{c}(n)})$}

For a topological space $X$ the quotient ring $H^{\ast }(X)/\tau (X)$ can be
identified with \textsl{the free part} of the integral cohomology $H^{\ast
}(X)$, and will be denoted by $\mathcal{F}(X)$. If $X=B_{SO(n)}$ and in
accordance to $n\equiv 0$ or $1$ we get from (2.3) that

\begin{quote}
$\mathcal{F}(B_{SO(n)})=\mathbb{Z}[p_{1},p_{2},\cdots ,p_{\left[ \frac{n-1}{2%
}\right] },e_{n}]$ or $\mathbb{Z}[p_{1},p_{2},\cdots ,p_{\left[ \frac{n-1}{2}%
\right] }]$,
\end{quote}

\noindent on which the map $\pi ^{\ast }$ injects by Lemma 4.1.

Consider the graded free $\mathbb{Z}$--module over the ring $\pi ^{\ast }%
\mathcal{F}(B_{SO(n)})$

\begin{enumerate}
\item[(5.1)] $A:=\pi ^{\ast }\mathcal{F}(B_{SO(n)})\otimes \Delta
(q_{0},\cdots ,q_{h(n)-1})\otimes \mathbb{Z}[\theta _{n}]$,
\end{enumerate}

\noindent and let $\varphi :A\rightarrow H^{\ast }(B_{Spin^{c}(n)})$ be the
additive homomorphism induced by the inclusions $\theta _{n},q_{r}\in
H^{\ast }(B_{Spin^{c}(n)})$.

\bigskip

\noindent \textbf{Lemma 5.1.} \textsl{The map }$\varphi $ \textsl{carries }$%
A $\textsl{\ isomorphically onto }$\mathcal{F}(B_{Spin^{c}(n)})$\textsl{.}

\bigskip

\noindent \textbf{Proof. }For a prime $p$ denote by $\varphi _{p}$ the
induced map

\begin{quote}
$\varphi _{p}:=\varphi \otimes id:A\otimes \mathbb{Z}_{p}\rightarrow H^{\ast
}(B_{Spin^{c}(n)})\otimes \mathbb{Z}_{p}$.
\end{quote}

\noindent in which we already have, by $2\tau (B_{Spin^{c}(n)})=0$, that

\begin{enumerate}
\item[(5.2)] $H^{\ast }(B_{Spin^{c}(n)})\otimes \mathbb{Z}_{2}=H_{\beta
}^{\ast }(B_{Spin^{c}(n)})\oplus \tau (B_{Spin^{c}(n)})$,

\item[(5.3)] $H^{\ast }(B_{Spin^{c}(n)})\otimes \mathbb{Z}_{p}=H^{\ast
}(B_{Spin^{c}(n)};\mathbb{Z}_{p})$\textsl{\ }if\textsl{\ }$p\neq 2$.
\end{enumerate}

\noindent It suffices to show that $\varphi _{p}$ injects for $p=2$, and is
an isomorphism for $p\neq 2$.

If $p=2$ we get from $\pi ^{\ast }\mathcal{F}(B_{SO(n)})\otimes \mathbb{Z}%
_{2}=\pi ^{\ast }H_{\beta }^{\ast }(B_{SO(n)})$ and (4.10) that

\begin{quote}
$A\otimes \mathbb{Z}_{2}=H_{\beta }^{\ast }(B_{Spin^{c}(n)})$.
\end{quote}

\noindent Comparing this with (5.2) shows that $\varphi _{2}$ injects.

Assume next that\textbf{\ }$p\neq 2$. With $2\tau (B_{SO(n)})=0$ we have $%
\pi ^{\ast }\mathcal{F}(B_{SO(n)})\otimes \mathbb{Z}_{p}=\pi ^{\ast }H^{\ast
}(B_{SO(n)};\mathbb{Z}_{p})$ and consequently

\begin{quote}
$A\otimes \mathbb{Z}_{p}=\pi ^{\ast }H^{\ast }(B_{SO(n)};\mathbb{Z}%
_{p})\otimes \Delta _{p}(q_{0},\cdots ,q_{h(n)-1})\otimes \mathbb{Z}%
_{p}[\theta _{n}]$.
\end{quote}

\noindent Since the composition $i^{\ast }\circ \psi _{p}$ carries the
factor $\Delta _{p}(q_{0},\cdots ,q_{h(n)-1})\otimes \mathbb{Z}_{p}[\theta
_{n}]$ isomorphically onto $H^{\ast }(\mathbb{C}P^{\infty };\mathbb{Z}_{p})$
by Lemma 4.4, we obtain

\begin{quote}
$A\otimes \mathbb{Z}_{p}=\pi ^{\ast }H^{\ast }(B_{SO(n)},\mathbb{Z}%
_{p})\otimes H^{\ast }(\mathbb{C}P^{\infty };\mathbb{Z}_{p})=H^{\ast
}(B_{Spin^{c}(n)};\mathbb{Z}_{p})$,
\end{quote}

\noindent where the second equality follows from Lemma 4.1. Comparing this
with (5.3) shows that $\psi _{p}$ is an isomorphism for $p\neq 2$.$\square $

\bigskip

Combining Lemma 5.1 with formula (4.4) we can express $H^{\ast
}(B_{Spin^{c}(n)})$ in term of its free part $\mathcal{F}(B_{Spin^{c}(n)})$
and the torsion ideal $\tau (B_{Spin^{c}(n)})$.

\bigskip

\noindent \textbf{Theorem 5.2. }\textsl{The cohomology }$H^{\ast
}(B_{Spin^{c}(n)})$\textsl{\ has the additive presentation}

\begin{enumerate}
\item[(5.4)] $\pi ^{\ast }\mathcal{F}(B_{SO(n)})\otimes \Delta (q_{0},\cdots
,q_{h(n)-1})\otimes \mathbb{Z}[\theta _{n}]\oplus \pi ^{\ast }\tau
(B_{SO(n)})\otimes \mathbb{Z}[\theta _{n}]$\textsl{,}
\end{enumerate}

\noindent \textsl{where the cup products between the free part and torsion
ideal are given by}

\begin{quote}
$\pi ^{\ast }(p_{i})\cup \pi ^{\ast }(z)=\pi ^{\ast }(w_{2i}^{2}\cup z)$ 
\textsl{(by }$\rho _{2}(p_{i})=w_{2i}^{2}$\textsl{) }

$q_{r}\cup \pi ^{\ast }(z)=\pi ^{\ast }(z\cup w_{2}^{(r)})$ \textsl{(by }$%
\rho _{2}(q_{r})=w_{2}^{(r)}$\textsl{),}
\end{quote}

\noindent \textsl{and} \textsl{where }$z\in \tau (B_{SO(n)})$, $1\leq i\leq %
\left[ \frac{n-1}{2}\right] $, $0\leq r\leq h(n)-1$.$\square $

\bigskip

By virtue of Theorem 5.2 we clarify a crucial relation relating the class $%
q_{h(n)-1}$ with the Euler class $\theta _{n}$.

\bigskip

\noindent \textbf{Lemma 5.3. }\textsl{There exists an element}

\begin{quote}
$\alpha _{n}\in \pi ^{\ast }H^{\ast }(B_{SO(n)})\otimes \Delta (q_{0},\cdots
,q_{h(n)-1})$\textsl{, }$\deg \alpha _{n}=2^{h(n)+1}$\textsl{,\ }
\end{quote}

\noindent \textsl{so that the following relation\ holds on the ring }$%
H^{\ast }(B_{Spin^{c}(n)})$

\begin{enumerate}
\item[(5.5)] $4(-1)^{h(n)}\theta _{n}+q_{h(n)-1}^{2}=a_{n}$\textsl{.}
\end{enumerate}

\noindent \textsl{In particular, the class }$a_{n}$\textsl{\ is
characterized uniquely by}

\begin{quote}
$\rho (a_{n})=(4(-1)^{h(n)}\rho _{0}(\theta
_{n})+g_{h(n)-1}^{2},(w_{2}^{(h(n)-1)})^{2})$\textsl{\ (see Lemma 4.1).}
\end{quote}

\noindent \textbf{Proof. }Since in degree $2^{h(n)+1}$ the quotient of $A$
by its subgroup $\pi ^{\ast }\mathcal{F}(B_{SO(n)})\otimes \Delta
(q_{0},\cdots ,q_{h(n)-1})\otimes 1$ is isomorphic to $\mathbb{Z}$ with the
basis $\{\theta _{n}\}$, formula (5.4) implies that the integral class $%
q_{h(n)}\in H^{\ast }(B_{Spin^{c}(n)})$ admits the expression

\begin{quote}
$q_{h(n)}=k\cdot \theta _{n}+\varphi (\beta _{n})$ for some $\beta _{n}\in
\pi ^{\ast }\mathcal{F}(B_{SO(n)})\otimes \Delta (q_{0},\cdots ,q_{h(n)-1})$
\end{quote}

\noindent and some integer $k\in \mathbb{Z}$. Applying $i^{\ast }$ we get by
Lemma 4.4 that

\begin{quote}
$2(-1)^{h(n)}x^{2^{h(n)}}=k\cdot x^{2^{h(n)}}$ on $H^{2^{h(n)+1}}(\mathbb{C}%
P^{\infty })$,
\end{quote}

\noindent showing $k=2(-1)^{h(n)}$. On the other hand, by iii) of Theorem C
we have

\begin{quote}
$2q_{h(n)}+q_{h(n)-1}^{2}=\pi ^{\ast }f(w_{2}^{(h(n)-1)})\in \pi ^{\ast
}H^{+}(B_{SO(n)})$.
\end{quote}

\noindent We obtain (5.5) by taking $a_{n}=\pi ^{\ast
}f(w_{2}^{(h(n)-1)})\otimes 1-2\beta _{n}$.$\square $

\bigskip

\noindent \textbf{Proof of Theorem D.} Set $\widetilde{A}:=\pi ^{\ast
}H^{\ast }(B_{SO(n)})\otimes \mathbb{Z}[q_{0},q_{1},\cdots
,q_{h(n)-1},\theta _{n}]$, and let $\widetilde{\varphi }:\widetilde{A}%
\rightarrow H^{\ast }(B_{Spin^{c}(n)})$ be induced by the inclusions $%
q_{r},\theta _{n}\in H^{\ast }(B_{Spin^{c}(n)})$. Then by the properties ii)
and iii) of Theorem C the map $\widetilde{\varphi }$ satisfies that

\begin{enumerate}
\item[(5.6)] $\widetilde{\varphi }(2\otimes q_{r+1}+1\otimes q_{r}^{2}-\pi
^{\ast }f(w_{2}^{(r)})\otimes 1)=0$, $0\leq r\leq h(n)-2$,

\item[(5.7)] $\widetilde{\varphi }(\pi ^{\ast }\delta _{2}(x)\otimes
q_{r}-\pi ^{\ast }\delta _{2}(x\cup w_{2}^{(r)})\otimes 1)=0$.
\end{enumerate}

\noindent In addition, we have by the relation (5.5) that

\begin{enumerate}
\item[(5.8)] $\widetilde{\varphi }(4(-1)^{h(n)}\otimes \theta _{n}+1\otimes
q_{h(n)-1}^{2}-a_{n})=0$.
\end{enumerate}

\noindent Therefore $R_{n}\subseteq \ker \widetilde{\varphi }$, where $R_{n}$
is the ideal on $\widetilde{A}$ stated in Theorem D. It remains to show that
the map $\widetilde{\varphi }$, when factoring through the quotient $%
\widetilde{A}\func{mod}R_{n}$, yields an additive isomorphism onto $H^{\ast
}(B_{Spin^{c}(n)})$.

Putting in transparent forms the relations (5.6) and (5.8) tell that

\begin{quote}
$q_{r}^{2}=-2q_{r+1}+\pi ^{\ast }f(w_{2}^{(r)})$, $0\leq r\leq h(n)-2;$

$q_{h(n)-1}^{2}=4(-1)^{h(n)+1}\theta _{n}+a_{n}$,
\end{quote}

\noindent respectively. These are the relations required to reduce each
monomial

\begin{quote}
$q_{0}^{r_{0}}q_{1}^{r_{1}}\cdots q_{h(n)-1}^{r_{h(n)-1}}\in \widetilde{A}$
with $r_{i}\geq 2$ for some $i$
\end{quote}

\noindent to an element of the $\pi ^{\ast }H^{\ast }(B_{SO(n)})$ module

\begin{quote}
$A_{1}:=\pi ^{\ast }H^{\ast }(B_{SO(n)})\otimes \Delta (q_{0},\cdots
,q_{h(n)-1})\otimes \mathbb{Z}[\theta _{n}]$.
\end{quote}

\noindent Moreover, combining $H^{\ast }(B_{SO(n)})=\mathcal{F}%
(B_{SO(n)})\oplus \tau (B_{SO(n)})$ with Lemma 5.1 we get the decomposition

\begin{quote}
$A_{1}=\mathcal{F}(B_{Spin^{c}(n)})\oplus \pi ^{\ast }\tau
(B_{SO(n)})\otimes \Delta (q_{0},\cdots ,q_{h(n)-1})\otimes \mathbb{Z}%
[\theta _{n}]$.
\end{quote}

\noindent in which the relation (5.7) tells, for any $z\in H^{\ast
}(B_{SO(n)};\mathbb{Z}_{2})$, that

\begin{quote}
$\widetilde{\varphi }(\pi ^{\ast }\delta _{2}(z)\otimes q_{r})=\pi ^{\ast
}\delta _{2}(z)\cup q_{r}=\pi ^{\ast }\delta _{2}(z\cup w_{2}^{(r)})$ (by $%
\rho (q_{r})=w_{2}^{(r)}$),
\end{quote}

\noindent showing

\begin{quote}
$\widetilde{\varphi }(\pi ^{\ast }\tau (B_{SO(n)})\otimes \Delta
(q_{0},\cdots ,q_{h(n)-1})\otimes \mathbb{Z}[\theta _{n}])$ $=\tau
(B_{Spin^{c}(n)})$
\end{quote}

\noindent by (4.5). The proof of Theorem D has now been completed by Theorem
5.2.$\square $

\bigskip

\noindent \textbf{Example 5.4. }In concrete situations Theorem D is directly
applicable to present the ring $H^{\ast }(B_{Spin^{c}(n)})$ by a minimal
system of explicit generators and relations. For instance if $n=8$ we get by
Lemma 5.1 that

\begin{quote}
$\mathcal{F}(B_{Spin^{c}(8)})=\mathbb{Z}[p_{1},p_{2},p_{3},e_{8},\theta
_{8}]\otimes \Delta (q_{0},q_{1},q_{2})$

$=\mathbb{Z}[q_{0},q_{1},q_{2},p_{3},e_{8},\theta _{8}]/\left\langle 4\theta
_{8}-q_{2}^{2}-a_{8}\right\rangle $,
\end{quote}

\noindent where the second equality is obtained by subtracting the
generators $p_{1},p_{2}$ using the relations

\begin{quote}
$p_{1}=2q_{1}+q_{0}^{2}$, $p_{2}=2q_{2}+q_{1}^{2}$ by (5.6).
\end{quote}

\noindent Furthermore, by (4.5), as well as the obvious relation on $\tau
(B_{Spin(8)})$

\begin{center}
$q_{0}^{r_{1}}\cup q_{1}^{r_{2}}\cup q_{2}^{r_{3}}\cup p_{3}^{r_{4}}\cup
e_{8}^{r_{5}}\cup \delta _{2}(w_{6})=\delta _{2}(w_{2}^{r_{1}}\cup
w_{4}^{r_{2}}\cup w_{6}^{2r_{2}}\cup w_{8}^{r_{3}+r_{5}}\cup w_{6})$,
\end{center}

\noindent we find that $\tau (B_{Spin(8)})$ is the ideal on $H^{\ast
}(B_{Spin(8)})$ generated by the single integral class $\delta _{2}(w_{6})$
with degree $7$. Summarizing we get that

\begin{enumerate}
\item[(5.9)] $H^{\ast }(B_{Spin^{c}(8)})=\mathbb{Z}%
[q_{0},q_{1},q_{2},p_{3},e_{8},\theta _{8},\delta _{2}(w_{6})]/\left\langle
2\delta _{2}(w_{6}),4\theta _{8}-q_{2}^{2}-a_{8}\right\rangle $
\end{enumerate}

\noindent where computing with the invariant $\theta _{8}$ shows that

\begin{quote}
$a_{8}=e_{8}^{2}-2e_{8}q_{2}-q_{0}^{2}p_{3}+2e_{8}q_{0}^{2}q_{1}$.$\square $
\end{quote}

\section{The cohomology $H^{\ast }(B_{Spin(n)})$}

The method developed in the preceding sections \S 4--5 applies equally well
to formulate the ring $H^{\ast }(B_{Spin(n)})$. To clarify this consider the
fibration induced by the two sheets covering $Spin(n)\rightarrow SO(n)$

\begin{enumerate}
\item[(6.1)] $\mathbb{R}P^{\infty }\overset{i}{\rightarrow }B_{Spin(n)}%
\overset{\overline{\pi }}{\rightarrow }B_{SO(n)}$,
\end{enumerate}

\noindent where $\overline{\pi }=\pi \circ \psi $ in the notion of (2.5) and
(2.6).

\bigskip

\noindent \textbf{Theorem C'. }\textsl{There exists a unique sequence }$%
\left\{ \overline{q}_{r}\QTR{sl}{,\ }r\geq 1\right\} $ \textsl{of} \textsl{%
integral} \textsl{cohomology classes on }$B_{Spin(n)}$\textsl{, }$\deg 
\overline{q}_{r}=2^{r+1}$\textsl{,} \textsl{that} \textsl{satisfies the
following system}

\begin{quote}
\textsl{i)} $\rho _{2}(\overline{q}_{r})=\overline{\pi }^{\ast }w_{2}^{(r)}$%
\textsl{;}

\textsl{ii)} $2\overline{q}_{1}=\overline{\pi }^{\ast }p_{1}$\textsl{,} $2%
\overline{q}_{r+1}+\overline{q}_{r}^{2}=\overline{\pi }^{\ast
}f(w_{2}^{(r)}) $\textsl{, }$r\geq 1$\textsl{.}
\end{quote}

\noindent \textbf{Proof.} The existence of a sequence $\left\{ \overline{q}%
_{r}\QTR{sl}{,\ }r\geq 0\right\} $ fulfilling the properties i) and ii) is
obtained from ii) and iii) of Theorem C by setting $\overline{q}_{r}:=$ $%
\psi ^{\ast }(q_{r})$. Note that $\overline{q}_{0}=0$, because $q_{0}$
belongs to the image of the transgression in the Serre spectral sequence of
the $U(1)$ fibration (2.6).

For the uniqueness of the sequence $\left\{ \overline{q}_{r},\QTR{sl}{\ }%
r\geq 1\right\} $ so obtained we make use of the map $\rho $ in Lemma 4.2
for the space $B_{Spin(n)}$. With $2\tau (B_{Spin(n)})=0$ it injects by \cite%
[30.6.]{BH}. On the other hand, the relations i) and ii) suffice to express
the images $\rho (\overline{q}_{r})$, $r\geq 1$. This competes the proof.$%
\square $

\bigskip

Let $\eta _{n}$ be the oriented real vector bundle on $B_{Spin(n)}$
associated to the real spin representation of the group $Spin(n)$, $\dim _{%
\mathbb{R}}\eta _{n}=$ $2^{k(n)+1}$\cite{ABS}, and denote by $\overline{%
\theta }_{n}\in H^{\ast }(B_{Spin(n)})$ its Euler class, where

\begin{quote}
$k(n)=\left[ \frac{n-1}{2}\right] $ if $n\func{mod}8\neq 0,1,7$, or $\left[ 
\frac{n-1}{2}\right] -1$ otherwise.
\end{quote}

\noindent It is known that (see \cite[Theorem 8.3]{BW} or \cite[Theorem 6.5]%
{Q}):

\bigskip

\noindent \textbf{Lemma 6.1.} \textsl{The maps }$\psi $ \textsl{in (2.6) and 
}$i$ \textsl{in (6.1)} \textsl{satisfy, respectively, that}

\begin{enumerate}
\item[(6.2)] $\psi ^{\ast }(\theta _{n})=\overline{\theta }_{n}$ \textsl{or }%
$\overline{\theta }_{n}^{2}$\textsl{, in accordance to }$n\func{mod}8\neq
0,1,7$,\textsl{\ or otherwise.}

\item[(6.3)] $i^{\ast }(\overline{\theta }_{n})=z^{2^{k(n)+1}}$\textsl{,} 
\textsl{where }$z$\textsl{\ is the generator of} $H^{1}(\mathbb{R}P^{\infty
};\mathbb{Z}_{2})=\mathbb{Z}_{2}$.$\square $
\end{enumerate}

Let $\{\overline{q}_{r},r\geq 1\}$ be the unique sequence\ obtained in
Theorem C', and let $\overline{\sigma }$ be the transgression in the
fibration (6.1). The counterparts of the Lemmas 4.2, 4.5 are summarized in
the following result.

\bigskip

\noindent \textbf{Lemma 6.2. }\textsl{For }$X=B_{Spin(n)}$\textsl{\ the
algebras} $H^{\ast }(X;\mathbb{Z}_{2})$, $H_{\beta }^{\ast }(X)$ \textsl{and 
}$\tau (X)$\textsl{\ have the following presentations}

\begin{enumerate}
\item[(6.4)] $H^{\ast }(B_{Spin(n)};\mathbb{Z}_{2})=\overline{\pi }^{\ast
}H^{\ast }(B_{SO(n)};\mathbb{Z}_{2})\otimes \mathbb{Z}_{2}[\overline{\theta }%
_{n}]$

\item[(6.5)] $H_{\beta }^{\ast }(B_{Spin(n)})=\overline{\pi }^{\ast
}H_{\beta }^{\ast }(B_{SO(n)})\otimes \Delta _{2}(\overline{q}_{1},\cdots ,%
\overline{q}_{k(n)-1})\otimes \mathbb{Z}_{2}[\overline{\theta }_{n}]$\textsl{%
;}

\item[(6.6)] $\tau (B_{Spin(n)})=\overline{\pi }^{\ast }\tau
(B_{SO(n)})\otimes \mathbb{Z}_{2}[\overline{\theta }_{n}]$\textsl{,}
\end{enumerate}

\noindent \textsl{respectively}, \textsl{where the map }$\overline{\pi }$ 
\textsl{in (6.1)} \textsl{induces the isomorphisms}

\begin{quote}
$\overline{\pi }^{\ast }H^{\ast }(B_{SO(n)};\mathbb{Z}_{2})=H^{\ast
}(B_{SO(n)};\mathbb{Z}_{2})/\left\langle \overline{\sigma }(z),\overline{%
\sigma }(z^{2}),\cdots ,\overline{\sigma }(z^{2^{k(n)-1}})\right\rangle $%
\textsl{;}

$\overline{\pi }^{\ast }H_{\beta }^{\ast }(B_{SO(n)})=H_{\beta }^{\ast
}(B_{SO(n)})/(w_{2}^{2})$\textsl{.}$\square $
\end{quote}

\bigskip

\noindent \textbf{Remark 6.3. }Formula\textbf{\ }(6.4) is due to Quillen 
\cite[Theorem 6.5]{Q}

In\textbf{\ }\cite[Proposition 6.1]{ABP} Anderson, Brown, and Peterson
obtained a formula for the Bockstein $H_{\beta }^{\ast }(B_{Spin(n)})$ in
the stable range $n=\infty $. The unstable instances has been calculated by
Kono \cite[\S 3]{Ko}, Benson and Wood \cite[\S 9]{BW}. Taking the advantage
of Theorem C' our formula (6.5) makes the structure of $H_{\beta }^{\ast
}(B_{Spin(n)})$ transparent, whose proof is similar to that of (4.10).$%
\square $

\bigskip

To formulate an additive presentation of the integral cohomology $H^{\ast
}(B_{Spin(n)})$, and to take the care of the overlap $\overline{q}_{1}=2%
\overline{\pi }^{\ast }p_{1}$ by ii) of Theorem C', we introduce the graded
quotient group

\begin{quote}
$A:=\overline{\pi }^{\ast }\mathcal{F}(B_{SO(n)})\otimes \Delta (\overline{q}%
_{1},\cdots ,\overline{q}_{k(n)-1})\otimes \mathbb{Z}[\overline{\theta }_{n}]%
\func{mod}(\overline{q}_{1}-2\overline{\pi }^{\ast }p_{1})$,
\end{quote}

\noindent and let $\varphi :A\rightarrow H^{\ast }(B_{Spin(n)})$ the
additive map induced by the inclusions $\overline{\theta }_{n},\overline{q}%
_{r}\in H^{\ast }(B_{Spin(n)})$. The same argument used in the proof of
Lemma 5.1 is applicable to show that $\varphi $ is an isomorphism onto the
free part $\mathcal{F}(B_{Spin(n)})$ of the ring $H^{\ast }(B_{Spin(n)})$.
Combining this with formula (6.6) we obtain the following result, which is
compatible with Theorem 5.2.

\bigskip

\noindent \textbf{Theorem 6.4.} \textsl{The cohomology }$H^{\ast
}(B_{Spin(n)})$\textsl{\ has the additive presentation}

\begin{enumerate}
\item[(6.7)] $A\oplus \overline{\pi }^{\ast }\tau (B_{SO(n)})\otimes \mathbb{%
Z}[\overline{\theta }_{n}]$\textsl{,}
\end{enumerate}

\noindent \textsl{where the cup products between the free part and torsion
ideal are given by}

\begin{quote}
$\overline{\pi }^{\ast }(p_{i})\cup \overline{\pi }^{\ast }(z)=\overline{\pi 
}^{\ast }(w_{2i}^{2}\cup z)$\textsl{\ (by }$\rho _{2}(p_{i})=$\textsl{\quad }%
$w_{2i}^{2}$\textsl{)}

$\overline{q}_{r}\cup \overline{\pi }^{\ast }(z)=\overline{\pi }^{\ast
}(z\cup w_{2}^{(r)})$ \textsl{(by }$\rho _{2}(\overline{q}_{r})=$\textsl{%
\quad }$\overline{\pi }^{\ast }(w_{2}^{(r)})$\textsl{),}
\end{quote}

\noindent \textsl{and} \textsl{where }$z\in \tau (B_{SO(n)})$, $1\leq i\leq %
\left[ \frac{n-1}{2}\right] $, $0\leq r\leq k(n)-1$.$\square $

\bigskip

Consider the integral class $\overline{q}_{k(n)}\in H^{\ast }(B_{Spin(n)})$.
Since in degree $2^{k(n)+1}$ the quotient of $A$ by its subgroup

\begin{quote}
$A_{1}:=\overline{\pi }^{\ast }\mathcal{F}(B_{SO(n)})\otimes \Delta (%
\overline{q}_{1},\cdots ,\overline{q}_{k(n)-1})\func{mod}(\overline{q}_{1}-2%
\overline{\pi }^{\ast }p_{1})$
\end{quote}

\noindent is isomorphic to $\mathbb{Z}$ with the basis $\{\overline{\theta }%
_{n}\}$, we must have the expression

\begin{enumerate}
\item[(6.8)] $\overline{q}_{k(n)}=\kappa \cdot \overline{\theta }%
_{n}+\varphi (b_{n})$ for some $\kappa \in \mathbb{Z}$ and $b_{n}\in A_{1}$
\end{enumerate}

\noindent by Theorem 6.4. However, in contrast to the proof of Lemma 5.3 the
fiber inclusion $i$ in (6.1) fails to be help to evaluate $\kappa $.
Nevertheless, resorting to the canonical subgroup $U(k)\subset Spin(2k)$ we
shall show in Lemma 7.10 that $\kappa =2(-1)^{k(n)-1}$. Thus, combining
(6.8) with the relation by Theorem C'

\begin{quote}
$2\overline{q}_{k(n)}+\overline{q}_{k(n)-1}^{2}=\overline{\pi }^{\ast
}f(w_{2}^{(k(n)-1)})\in $ $\overline{\pi }^{\ast }H^{+}(B_{SO(n)})$
\end{quote}

\noindent yields the following relation relating the class $\overline{q}%
_{k(n)-1}$ with the Euler class $\overline{\theta }_{n}$

\begin{enumerate}
\item[(6.9)] $4(-1)^{k(n)-1}\overline{\theta }_{n}+\overline{q}%
_{k(n)-1}^{2}=\varphi (\overline{b}_{n})$, $\overline{b}_{n}\in \overline{%
\pi }^{\ast }\overline{H}^{+}(B_{SO(n)})\otimes \Delta (\overline{q}%
_{1},\cdots ,\overline{q}_{k(n)-1})$.
\end{enumerate}

Granted with Theorem C', Theorem 6.4 and formula (6.9), the same argument
used in the proof of Theorem D is applicable to obtain the following result,
whose proof can therefore be omitted.

\bigskip

\noindent \textbf{Theorem D'. }\textsl{The ring }$H^{\ast }(B_{Spin(n)})$%
\textsl{\ has the presentation}

\begin{enumerate}
\item[(6.10)] $H^{\ast }(B_{Spin(n)})=\overline{\pi }^{\ast }H^{\ast
}(B_{SO(n)})\otimes \mathbb{Z}[\overline{q}_{1},\cdots ,\overline{q}%
_{k(n)-1},\overline{\theta }_{n}]/K_{n}$\textsl{,}
\end{enumerate}

\noindent \textsl{in which }$K_{n}$ \textsl{denotes the ideal generated by
the following elements}

\begin{quote}
\textsl{i)} $2\overline{q}_{1}-\overline{\pi }^{\ast }p_{1}$\textsl{,} $2%
\overline{q}_{r+1}+\overline{q}_{r}^{2}-\overline{\pi }^{\ast
}f(w_{2}^{(r)}) $\textsl{,\quad }$1\leq r\leq k(n)-2$\textsl{;}

\textsl{ii)} $\overline{\pi }^{\ast }\delta _{2}(y)\cup \overline{q}_{r}-%
\overline{\pi }^{\ast }\delta _{2}(y\cup w_{2}^{(r)}),$ $1\leq r\leq k(n)-2$%
\textsl{;}

\textsl{iii)} $4(-1)^{k(n)-1}\overline{\theta }_{n}+\overline{q}%
_{k(n)-1}^{2}-\overline{b}_{n}$\textsl{,}
\end{quote}

\noindent \textsl{where }$y\in H^{\ast }(B_{SO(n)},\mathbb{Z}_{2})$\textsl{,}
$\overline{b}_{n}\in \overline{\pi }^{\ast }\overline{H}^{+}(B_{SO(n)})%
\otimes \Delta (\overline{q}_{1},\cdots ,\overline{q}_{k(n)-1})$\textsl{.}$%
\square $

\bigskip

\noindent \textbf{Example 6.5. }If $n=8$ we get from formula (6.7) that

\begin{quote}
$\mathcal{F}(B_{Spin(8)})=\mathbb{Z}[p_{1},p_{2},p_{3},e_{8},\overline{%
\theta }_{8}]\otimes \Delta (\overline{q}_{1})\func{mod}(\overline{q}_{1}-2%
\overline{\pi }^{\ast }p_{1})$

$\qquad =\mathbb{Z}[\overline{q}_{1},p_{3},e_{8},\overline{\theta }_{8}]$
\end{quote}

\noindent where the second equality follows from

\begin{quote}
$p_{1}=2\overline{q}_{1}$ (by i) of Theorem D'),

$p_{2}=4\overline{\theta }_{8}+2e_{8}-\overline{q}_{1}^{2}$ (by iii) of
Theorem D').
\end{quote}

\noindent Furthermore, by formula (6.6), together with the relation on $%
H^{\ast }(B_{Spin(8)})$

\begin{quote}
$\overline{q}_{1}^{r_{1}}\cup p_{3}^{r_{2}}\cup e_{8}^{r_{3}}\cup \delta
_{2}(w_{6})=\delta _{2}(w_{4}^{r_{1}}\cup w_{6}^{2r_{2}}\cup
w_{8}^{r_{3}}\cup w_{6})$,
\end{quote}

\noindent the ideal $\tau (B_{Spin(8)})$ is generated by $\delta _{2}(w_{6})$%
. Summarizing we get

\begin{enumerate}
\item[(6.12)] $H^{\ast }(B_{Spin(8)})=\mathbb{Z}[\overline{q}%
_{1},p_{3},e_{8},\overline{\theta }_{8},\delta _{2}(w_{6})]/\left\langle
2\delta _{2}(w_{6})\right\rangle $.
\end{enumerate}

Similarly, one can show that

\begin{enumerate}
\item[(6.13)] $H^{\ast }(B_{Spin(12)})=\frac{\mathbb{Z}[\overline{q}_{1},%
\overline{q}_{2},\overline{q}_{3},\overline{q}_{4},p_{3},p_{5},e_{12},%
\overline{\theta }_{12},\delta _{2}(w_{6}),\delta _{2}(w_{10})]}{%
\left\langle 2\delta _{2}(w_{6}),2\delta
_{2}(w_{10}),R_{32},R_{64}\right\rangle },$
\end{enumerate}

\noindent where

\begin{quote}
$R_{32}=2\overline{q}_{4}+\overline{q}_{3}^{2}-\overline{\pi }^{\ast
}f(w_{2}^{(3)});$

$R_{64}=4\overline{\theta }_{12}+\overline{q}_{4}^{2}-\overline{b}_{8}$,
\end{quote}

\noindent and where the deficiencies $\overline{\pi }^{\ast }f(w_{2}^{(3)})$
and $\overline{b}_{8}$ can be computed explicitly (they are too lengthy to
be presented).$\square $

\section{The subgroup $U(k)\subset Spin(2k)$}

The unitary group $U(k)$ of order $k$ is contained in the special orthogonal
group $SO(2k)$ as the centralizer of the circle subgroup $diag\{z,\cdots
,z\} $ with $k$ copies of $z\in SO(2)$. It is also contained in the group $%
Spin(2k)$ as the centralizer of the circle subgroup

\begin{center}
$U(1)=\{\underset{1\leq i\leq k}{\Pi }(\cos t+\sin t\cdot e_{2i-1}e_{2i})\in
Spin(2k),t\in \mathbb{R}\}$ (see \cite[p.173]{BD}, \cite{DL}).
\end{center}

\noindent Let $\lambda _{0}:U(k)\rightarrow SO(2k)$ and $\lambda
:U(k)\rightarrow Spin(2k)$ be the inclusions, respectively, and define $%
\lambda ^{c}$ to be the homomorphism

\begin{quote}
$\lambda ^{c}:=\lambda \times _{\mathbb{Z}_{2}}id:U^{c}(k):=U(k)\times _{%
\mathbb{Z}_{2}}U(1)\rightarrow Spin^{c}(2k)$.
\end{quote}

\noindent These homomorphisms then fit into the commutative diagram

\begin{enumerate}
\item[(7.1)] $%
\begin{array}{ccccc}
U(k) & \overset{a^{\prime }}{\rightarrow } & U^{c}(k) & \overset{b^{\prime }}%
{\rightarrow } & U(k) \\ 
\lambda \downarrow &  & \lambda ^{c}\downarrow &  & \lambda _{0}\downarrow
\\ 
Spin(2k) & \overset{a}{\rightarrow } & Spin^{c}(2k) & \overset{b}{%
\rightarrow } & SO(2k)%
\end{array}%
$,
\end{enumerate}

\noindent where $a$ (resp. $a^{\prime }$) denotes the inclusion onto the
kernel of the epimorphism

\begin{quote}
$Spin(2k)\times _{\mathbb{Z}_{2}}U(1)\rightarrow U(1)$ (resp. $U(k)\times _{%
\mathbb{Z}_{2}}U(1)\rightarrow U(1)$),
\end{quote}

\noindent and where $b$ (resp. $b^{\prime }$) is the obvious epimorphism

\begin{quote}
$Spin(2k)\times _{\mathbb{Z}_{2}}U(1)\rightarrow SO(2k)$ (resp. $U(k)\times
_{\mathbb{Z}_{2}}U(1)\rightarrow U(k)$).
\end{quote}

\noindent It is easy to see that

\bigskip

\noindent \textbf{Lemma 7.1.} \textsl{The composition }$b\circ a$\textsl{\
is the universal covering on }$SO(2k)$\textsl{.}

\textsl{The composition }$b^{\prime }\circ a^{\prime }$\textsl{\ is the two
sheets covering of }$U(k)$\textsl{\ by its self with associated }$\func{mod}%
2 $\textsl{\ Euler class the generator of }$H^{1}(U(k);\mathbb{Z}_{2})=%
\mathbb{Z}_{2}$\textsl{.}$\square $

\bigskip

The homomorphisms $\lambda $ and $\lambda ^{c}$ are particularly relevant to
the understanding of the relationship between spinors and complex structures
on vector bundles \cite{ABS}. The purpose of this section is to determine
their induced maps on the cohomologies of the relevant classifying spaces.
The results will bring us direct recurrences to produce the integral Weyl
invariants of the groups $Spin^{c}(n)$ and $Spin(n)$, by the algorithms
illustrated in Examples 7.6 and 7.9.

\bigskip

\textbf{7.1. The conventions.} For a homomorphism $h:H\rightarrow G$ of Lie
groups denote by $B_{h}:B_{H}\rightarrow B_{G}$ its induced map on the
classifying spaces. The diagram (7.1) then gives rise to the commutative
diagram of induced ring maps

\begin{enumerate}
\item[(7.2)] $%
\begin{array}{ccccc}
H^{\ast }(B_{U(k)}) & \overset{B_{a^{\prime }}^{\ast }}{\leftarrow } & 
H^{\ast }(B_{U^{c}(k)}) & \overset{B_{b^{\prime }}^{\ast }}{\leftarrow } & 
H^{\ast }(B_{U(k)}) \\ 
B_{\lambda }^{\ast }\uparrow &  & B_{\lambda ^{c}}^{\ast }\uparrow &  & 
B_{\lambda _{0}}^{\ast }\uparrow \\ 
H^{\ast }(B_{Spin(2k)}) & \overset{\psi ^{\ast }}{\longleftarrow } & H^{\ast
}(B_{Spin^{c}(2k)}) & \overset{\pi ^{\ast }}{\longleftarrow } & H^{\ast
}(B_{SO(2k)})%
\end{array}%
$
\end{enumerate}

\noindent where $\psi =B_{a}$, $\pi =B_{b}$ in the notation of (2.5) and
(2.6). By the naturality of characteristic classes, and for the sake to
simplify notation, the following convention is adopted throughout this
section.

\bigskip

\textbf{Convention a).} Let $1+c_{1}+\cdots +c_{k}$ be the total Chern class
of the canonical bundle on $B_{U(k)}$, and recall that $H^{\ast }(B_{U(k)})=%
\mathbb{Z}[c_{1},\cdots ,c_{k}]$. The notion $c_{r}\in H^{\ast
}(B_{U^{c}(k)})$ will be used to simplify $B_{b^{\prime }}^{\ast }(c_{r})$.

\bigskip

\textbf{Convention b). }Let $p_{r}\in H^{\ast }(B_{SO(2k)})$ be the $r^{th}$
Pontryagin class if $r<k$, or $e_{2k}^{2}$ if $r=k$. The same notion $p_{r}$
will also be reserved for both of the classes

\begin{quote}
$\pi ^{\ast }p_{r}\in H^{\ast }(B_{Spin^{c}(2k)})$ and $\psi ^{\ast }\pi
^{\ast }p_{r}\in H^{\ast }(B_{Spin(2k)})$,
\end{quote}

\noindent as well as the polynomial $B_{\lambda _{0}}^{\ast }(p_{r})\in $ $%
H^{\ast }(B_{U(k)})$ and $B_{\lambda ^{c}}^{\ast }(p_{r})\in $ $H^{\ast
}(B_{U^{c}(k)})$ whose canonical expression is (where $c_{i}=0$\ if $i$\ $>k$%
)

\begin{enumerate}
\item[(7.3)] $p_{r}=c_{r}^{2}-2c_{r-1}c_{r+1}+\cdots
+2(-1)^{r-1}c_{1}c_{2r-1}+2(-1)^{r}c_{2r}$ (\cite[p.177]{MS}).
\end{enumerate}

\textbf{Convention c).} The notion $e_{2k}$ for the Euler class of the
canonical bundle on $B_{SO(2k)}$ is also applied to denote either

\begin{quote}
$\pi ^{\ast }e_{2k}\in H^{\ast }(B_{Spin^{c}(2k)})$ or $\psi ^{\ast }\pi
^{\ast }e_{2k}\in H^{\ast }(B_{Spin(2k)})$.
\end{quote}

\noindent These convention will not cause any confusion, because in each
circumstance the cohomologies or the homomorphisms involved will be clearly
stated.

\bigskip

\textbf{7.2. The operators }$(\psi ,\delta )$\textbf{\ on the} \textbf{%
cohomology }$H^{\ast }(B_{U^{c}(k)})$\textbf{.} The group $U^{c}(k)$ has two
obvious $1$--dimensional unitary representations

\begin{quote}
$U(k)\times _{\mathbb{Z}_{2}}U(1)\rightarrow U(k)\overset{\det }{\rightarrow 
}U(1)$ and $U(k)\times _{\mathbb{Z}_{2}}U(1)\rightarrow U(1)$
\end{quote}

\noindent whose Euler classes are $c_{1}$ and $B_{\lambda ^{c}}^{\ast
}(q_{0})$, respectively. By the commutivity of the second diagram in (7.2)
we have

\begin{quote}
$\rho _{2}(B_{\lambda ^{c}}^{\ast }(q_{0})+c_{1})=B_{\lambda ^{c}}^{\ast
}\circ B_{b}^{\ast }(w_{2})+B_{b^{\prime }}^{\ast }\circ B_{\lambda
_{0}}^{\ast }(w_{2})=0$,
\end{quote}

\noindent implying that the sum $B_{\lambda ^{c}}^{\ast }(q_{0})+c_{1}\in
H^{2}(B_{U^{c}(k)})$ is divisible by $2$. This brings us the integral class

\begin{quote}
$y:=\frac{1}{2}(B_{\lambda ^{c}}^{\ast }(q_{0})+c_{1})\in
H^{2}(B_{U^{c}(k)}) $
\end{quote}

\noindent by which following result becomes transparent by Lemma 7.1.

\bigskip

\noindent \textbf{Lemma 7.2.} \textsl{We have }$H^{\ast }(B_{U^{c}(k)})=%
\mathbb{Z}[y,c_{1},\cdots ,c_{k}]$\textsl{\ so that}

\begin{quote}
\textsl{i)} $B_{b^{\prime }}^{\ast }(c_{r})=c_{r},1\leq r\leq k$;

\textsl{ii) }$B_{a^{\prime }}^{\ast }(z)=c_{1},2c_{1}$\textsl{\ or }$c_{r}$%
\textsl{\ for }$z=y,c_{1}$\textsl{\ or }$c_{r}$\textsl{\ with }$r\geq 2$%
\textsl{;}

\textsl{iii) }$B_{\lambda ^{c}}^{\ast }(q_{0})=2y-c_{1}$\textsl{.}$\square $
\end{quote}

For each ordered sequence $\lambda =(\lambda _{1},\cdots ,\lambda _{k})$ of $%
k$ non-negative integers define

\begin{quote}
$c_{\lambda }:=c_{1}^{\lambda _{1}}\cdots c_{k}^{\lambda _{k}}\in H^{\ast
}(B_{U^{c}(k)})$.
\end{quote}

\noindent Since all the monomials $y^{r}c_{\lambda }$, $r\geq 0$, form an
additive basis of the cohomology $H^{\ast }(B_{U^{c}(k)})$ every element $%
u\in H^{\ast }(B_{U^{c}(k)})$ has the unique expansion

\begin{quote}
$u=\underset{(r,\lambda )}{\Sigma }u_{(r,\lambda )}\cdot y^{r}c_{\lambda }$, 
$u_{(r,\lambda )}\in \mathbb{Z}$.
\end{quote}

\noindent We may then introduce the operator $\psi $ on the ring $H^{\ast
}(B_{U^{c}(k)})$ by

\begin{quote}
$\psi (u):=\underset{(r,\lambda )}{\Sigma }\rho _{2}(u_{(r,\lambda )})\cdot
y^{2r}p_{\lambda },$ $u\in H^{\ast }(B_{U^{c}(k)})$,
\end{quote}

\noindent where $p_{\lambda }=p_{1}^{\lambda _{1}}\cdots p_{k}^{\lambda
_{k}} $ with $p_{i}$ the polynomials in the Chern classes $c_{r}$ given by
(7.3). Alternatively,

\bigskip

\noindent \textbf{Corollary 7.3.} \textsl{The operator }$\psi $\textsl{\ is
characterized uniquely by the following algorithmic properties, where }$%
u,v\in H^{\ast }(B_{U^{c}(k)})$\textsl{, }$b\in \mathbb{Z}$\textsl{,}

\begin{quote}
\textsl{i) }$\psi (u+b\cdot y^{r}c_{\lambda })=\psi (u)+\rho _{2}(b)\cdot
y^{2r}p_{\lambda }$\textsl{\ if }$\rho _{2}(u_{(r,\lambda )})=0;$

\textsl{ii) }$\psi (u\cup v)=\psi (u)\cup \psi (v)$\textsl{.}$\square $
\end{quote}

It follows from $p_{\lambda }\equiv c_{\lambda }^{2}\func{mod}2$ by (7.3)
that the operator $\psi $ satisfies also the relation $\psi (u)\equiv u^{2}%
\func{mod}2$, which implies that there exists a unique operator $\delta $ on
the ring $H^{\ast }(B_{U^{c}(k)})$ that is related to $\psi $ by the formula

\begin{quote}
$\psi (u)=u^{2}+2\delta (u).$
\end{quote}

\noindent \textbf{Definition 7.4. }For a polynomial $u\in
H^{2n}(B_{U^{c}(k)})$ the sequence $\left\{ u,\delta (u),\delta
^{2}(u),\cdots \right\} $ on $H^{\ast }(B_{U^{c}(k)})$ defined inductively
by $\delta ^{r}(u)=\delta (\delta ^{r-1}(u)$ will be called \textsl{the
derived sequence} \textsl{of} \textsl{the initial polynomial} $u$.$\square $

\bigskip

Since the ring $H^{\ast }(B_{U^{c}(k)})$ is torsion free we get

\bigskip

\noindent \textbf{Corollary 7.5. }\textsl{For any }$u\in H^{\ast
}(B_{U^{c}(k)})$\textsl{\ the derive sequence }$\left\{ u,\delta (u),\delta
^{2}(u),\cdots \right\} $\textsl{\ can be computed by }$\psi $ \textsl{via
the following recurrence relations}

\begin{quote}
$\delta ^{0}(u)=u$\textsl{, }$2\delta ^{r+1}(u)+\delta ^{r}(u)^{2}=\psi
(\delta ^{r}(u))$, $r\geq 0$.$\square $
\end{quote}

\bigskip

\noindent \textbf{Example 7.6. }Corollaries 7.3 and 7.5 indicate a direct
and effective recurrence to produce the sequence $\left\{ u,\delta
(u),\delta ^{2}(u),\cdots \right\} $ from the initial one $u$. For example
we take $u=2y-c_{1}\in H^{2}(B_{U^{c}(k)})$. Then

\begin{quote}
$\delta ^{1}(u)=-c_{2}+2yc_{1}-2y^{2}$;

$\delta
^{2}(u)=c_{4}-c_{1}c_{3}-2y^{2}c_{1}^{2}-2y^{4}+2yc_{1}c_{2}-2y^{2}c_{2}+4y^{3}c_{1} 
$;

$\delta
^{3}(u)=c_{_{8}}-c_{_{1}}c_{_{7}}+c_{2}c_{6}-c_{3}c_{5}+(c_{1}^{2}-2c_{2})(-c_{2}c_{4}+c_{1}c_{5}-c_{6}) 
$

$\qquad
-c_{2}c_{3}^{2}-c_{1}c_{3}c_{4}-2(-y^{2}c_{1}^{2}-y^{4}+yc_{1}c_{2}-y^{2}c_{2}+2y^{3}c_{1})^{2} 
$

$\qquad
-2(c_{4}-c_{1}c_{3})(-y^{2}c_{1}^{2}-y^{4}+yc_{1}c_{2}-y^{2}c_{2}+2y^{3}c_{1}) 
$, $\cdots $.
\end{quote}

\noindent In addition, it can be shown that

\begin{enumerate}
\item[(7.4)] $\psi (\delta ^{r}(u))\in \left\langle p_{2},\cdots
,p_{k}\right\rangle $, $r>1$.$\square $
\end{enumerate}

\textbf{7.3. The ring map }$B_{\lambda ^{c}}^{\ast }:H^{\ast
}(B_{Spin^{c}(2k)})\rightarrow H^{\ast }(B_{U^{c}(k)})$\textbf{.} By the
conventions in \textbf{7.1} and with respect to the presentation

\begin{quote}
$H^{\ast }(B_{Spin^{c}(2k)})=\pi ^{\ast }H^{\ast }(B_{SO(2k)})\otimes 
\mathbb{Z}[q_{0},q_{1},\cdots ,q_{h(2k)-1},\theta _{2k}]/R_{2k}$,
\end{quote}

\noindent by Theorem D, partial information on the ring map $B_{\lambda
^{c}}^{\ast }$ has already known. Precisely we have

\begin{quote}
$B_{\lambda ^{c}}^{\ast }(e_{2k})=c_{k}$, \textsl{\ }$B_{\lambda ^{c}}^{\ast
}(p_{r})=p_{r}$, $B_{\lambda ^{c}}^{\ast }(\tau (B_{Spin^{c}(2k)}))=0$,
\end{quote}

\noindent where the third relation follows from the fact that the target
ring $H^{\ast }(B_{U^{c}(k)})$ is torsion free. In addition, since $\theta
_{n}\in H^{\ast }(B_{Spin^{c}(n)})$ is the Euler class of the the complex
spin presentation of the group $Spin^{c}(n)$, the polynomials $B_{\lambda
^{c}}^{\ast }(\theta _{n})$ can be computed in representation theory \cite%
{ABS}. As examples we have

\begin{quote}
$B_{\lambda ^{c}}^{\ast }(\theta _{4})=y^{2}-yc_{1};$

$B_{\lambda ^{c}}^{\ast }(\theta
_{6})=y^{4}-2y^{3}c_{1}+y^{2}(c_{2}+c_{1}^{2})-y(c_{1}c_{2}-c_{3})$;

$B_{\lambda ^{c}}^{\ast }(\theta
_{8})=y^{8}-4y^{7}c_{1}+y^{6}(2c_{2}+6c_{1}^{2})-y^{5}(6c_{3}+4c_{1}^{3}+6c_{1}c_{2}) 
$

$\
+y^{4}(c_{1}^{4}+c_{2}^{2}+6c_{1}^{2}c_{2}+c_{1}c_{3}-4c_{4})-y^{3}(2c_{1}c_{2}^{2}-8c_{1}c_{4}+2c_{1}^{2}c_{3}+2c_{1}^{3}c_{2}) 
$

\ \ $%
+y^{2}(c_{1}^{3}c_{3}-22c_{2}c_{4}+c_{1}c_{2}c_{3}-c_{3}^{2}-5c_{1}^{2}c_{4}+\allowbreak c_{1}^{2}c_{2}^{2})-y(c_{1}^{2}c_{2}c_{3}-c_{1}^{3}c_{4}-c_{1}c_{3}^{2}) 
$.
\end{quote}

\noindent Summarizing, the determination of the ring map $B_{\lambda
^{c}}^{\ast }$ is reduced to express the sequence $\{B_{\lambda ^{c}}^{\ast
}(q_{r}),0\leq r\leq h(n)-1\}$ as explicit polynomials in $H^{\ast
}(B_{U^{c}(k)})$.

\bigskip

\noindent \textbf{Theorem 7.7.} \textsl{The sequence }$\{B_{\lambda
^{c}}^{\ast }(q_{r}),0\leq r\leq h(n)-1\}$\textsl{\ is the first }$h(n)$ 
\textsl{terms of the} \textsl{derived sequence }$\left\{ u,\delta (u),\delta
^{2}(u),\cdots \right\} $ \textsl{of} $u=c_{1}-2y$\textsl{\ (see in Example
7.6).}

\bigskip

\noindent \textbf{Proof.} The proof serves also the purpose to bring a
passage from the operators $\{f,\gamma \}$ on $H^{\ast }(B_{SO(2k)};\mathbb{Z%
}_{2})$ (Theorem B) to the ones $\{\psi ,\delta \}$ on $H^{\ast
}(B_{U^{c}(k)})$ via the map $B_{\lambda ^{c}}^{\ast }$. In view of the
presentation

\begin{quote}
$H^{\ast }(B_{U^{c}(k)};\mathbb{Z}_{2})=\mathbb{Z}_{2}[\rho _{2}(y),\rho
_{2}(c_{1}),\cdots ,\rho _{2}(c_{k})]$
\end{quote}

\noindent by Lemma 7.2 define the operator $R:H^{2k}(B_{U^{c}(k)};\mathbb{Z}%
_{2})\rightarrow H^{4k}(B_{U^{c}(k)})$ by

\begin{enumerate}
\item[(7.5)] $R(u):=\left\{ 
\begin{tabular}{l}
$y^{2r}p_{\lambda }$ if $u=\rho _{2}(y^{r}c_{\lambda })$ is a monomial; \\ 
$R(u_{1})+\cdots +R(u_{m})$ if $u=$ $u_{1}+\cdots +u_{m}$%
\end{tabular}%
\right. $
\end{enumerate}

\noindent where the $u_{i}$'s are distinct monomials in the $\rho
_{2}(y),\rho _{2}(c_{1}),\cdots ,\rho _{2}(c_{k})$. Then, in addition to the
obvious factorization

\begin{enumerate}
\item[(7.6)] $\psi =R\circ \rho _{2}$,
\end{enumerate}

\noindent comparison with the formula (2.4) of the operator $f$ tells that

\begin{enumerate}
\item[(7.7)] $R\circ B_{\lambda ^{c}}^{\ast }\circ \pi ^{\ast }=B_{\lambda
^{c}}^{\ast }\circ \pi ^{\ast }\circ f$ (where $\pi ^{\ast }=B_{b}^{\ast }$).
\end{enumerate}

\noindent On the other hand, by the second diagram in (7.2), applying the
ring map $B_{\lambda ^{c}}^{\ast }$ to the relation ii) and iii) of Theorem
C gives rise to

\begin{enumerate}
\item[(7.8)] $\rho _{2}\circ B_{\lambda ^{c}}^{\ast }(q_{r})=B_{\lambda
^{c}}^{\ast }\circ \pi ^{\ast }(w_{2}^{(r)})$ on $H^{\ast }(B_{U^{c}(k)};%
\mathbb{Z}_{2})$, and

\item[(7.9)] $2B_{\lambda ^{c}}^{\ast }(q_{r+1})+B_{\lambda ^{c}}^{\ast
}(q_{r})^{2}=B_{\lambda ^{c}}^{\ast }\circ \pi ^{\ast }f(w_{2}^{(r)})$, $%
r\geq 0$,\textsl{\ }on $H^{\ast }(B_{U^{c}(k)})$,
\end{enumerate}

\noindent respectively. Setting $u_{r+1}:=B_{\lambda ^{c}}^{\ast }(q_{r})$
these imply that

\begin{quote}
$2u_{r+1}+u_{r}^{2}=B_{\lambda ^{c}}^{\ast }\circ \pi ^{\ast }\circ
f(w_{2}^{(r)})$ (by (7.9))

$=R\circ B_{\lambda ^{c}}^{\ast }\circ \pi ^{\ast }(w_{2}^{(r)})$ (by (7.7))

$=R\circ \rho _{2}(u_{r})$ (by (7.8))

$=\psi (u_{r})$ (by (7.6)).
\end{quote}

\noindent With $u_{1}=2y-c_{1}$ by iii) of Lemma 7.2 we obtain the result by
Corollary 7.5.$\square $

\bigskip

\textbf{7.4. The ring map }$B_{\lambda }^{\ast }:H^{\ast
}(B_{Spin(2k)})\rightarrow H^{\ast }(B_{U(k)})$\textbf{. }Carrying on
discussion from Example 7.6 let $\left\{ u,\delta (u),\delta ^{2}(u),\cdots
\right\} $ be the derived sequence of the polynomial $u=2y-c_{1}$. As
elements in the polynomial ring $H^{\ast }(B_{U^{c}(k)})$ we can write for
each $r\geq 1$ that

\begin{quote}
$\delta ^{r}(u)=u^{(r)}(y,c_{1},c_{2},\cdots ,c_{k})$.
\end{quote}

\noindent Applying the ring map $B_{a^{\prime }}^{\ast }$ to the polynomials 
$p_{r},\delta ^{r}(u)$, $\psi (\delta ^{r}(u))\in H^{\ast }(B_{U^{c}(k)})$
we obtain by ii) of Lemma 7.2 the following polynomials in $H^{\ast
}(B_{U(k)}):$

\begin{enumerate}
\item[(7.10)] $g_{r}:=B_{a^{\prime }}^{\ast
}(p_{r})=p_{r}+2(-1)^{r-1}c_{1}c_{2r-1}$ (see (7.3) for $p_{r}$);

\item[(7.11)] $\alpha _{r}:=B_{a^{\prime }}^{\ast }(\delta
^{r}(u))=u^{(r)}(c_{1},2c_{1},c_{2},\cdots ,c_{k})$,

\item[(7.12)] $f_{r}:=B_{a^{\prime }}^{\ast }(\psi (\delta
^{r}(u))=u^{(r)}(g_{1},2g_{1},g_{2},\cdots ,g_{k})$
\end{enumerate}

\noindent On the other hand, according to Theorem D' and by the conventions
in 7.1, the ring $H^{\ast }(B_{Spin(2k)})$ is generated multiplicatively by
the integral classes

\begin{quote}
$\{p_{r},$ $1\leq r\leq k-1\}$; $\{\overline{q}_{r},1\leq r\leq h(n)-1\}$; $%
e_{2k}$, $\overline{\theta }_{n}$,
\end{quote}

\noindent together with the ideal $\tau (B_{Spin^{c}(n)})$. Thus, we obtain
from Theorem 7.7, the relation $\overline{q}_{r}=\psi ^{\ast }(q_{r})$ by
Theorem C', as well as the commutivity of the first diagram in (7.2), the
following result.

\bigskip

\noindent \textbf{Theorem 7.8.} \textsl{The map }$B_{\lambda }^{\ast
}:H^{\ast }(B_{Spin(2k)})\rightarrow H^{\ast }(B_{U(k)})$\textsl{\ is
determined by}

\begin{quote}
\textsl{i)} $B_{\lambda }^{\ast }(p_{r})=g_{r},$ $1\leq r\leq k-1$;

\textsl{ii)} $B_{\lambda }^{\ast }(\overline{q}_{r})=\alpha _{r},1\leq r\leq
h(n)-1$,
\end{quote}

\noindent \textsl{and by }$B_{\lambda }^{\ast }(e_{2k})=c_{k},$ $B_{\lambda
}^{\ast }(\tau (B_{Spin^{c}(n)}))=0$\textsl{.}$\square $

\bigskip

\noindent \textbf{Example 7.9. }We shall show in Section \S 8 that the
sequence $\{\alpha _{r},1\leq r\leq h(n)-1\}$ (as well as $\{p_{r},1\leq
r\leq k-1\}$) of polynomials in the Chern classes are integral Weyl
invariants of the group $Spin(n)$. We emphasize at this stage that these
polynomials can be effectively produced by the simple algorithm indicated by
Corollary 7.5, together with the formula (7.11). As examples, combining
results in Example 7.6 with formula (7.11) we obtain that

\begin{quote}
$\alpha _{1}=-c_{2}+2c_{1}^{2}$;

$\alpha _{2}=c_{4}-2c_{1}c_{3}+2c_{1}^{2}c_{2}-2c_{1}^{4}$;

$\alpha
_{3}=-c_{8}+2c_{1}c_{7}-c_{2}c_{6}+c_{3}c_{5}-(2c_{1}^{2}-c_{2})(c_{3}^{2}-2c_{2}c_{4}+4c_{1}c_{5}-2c_{6}) 
$

$\qquad
-2c_{4}(c_{1}c_{3}-c_{1}^{2}c_{2}+c_{1}^{4})+2(c_{1}c_{3}-c_{1}^{2}c_{2}+c_{1}^{4})^{2},\cdots 
${\small .}$\square $
\end{quote}

\bigskip

\textbf{7.5.} We conclude this section with the following result which has
played a role in showing Theorem D'.

\bigskip

\noindent \textbf{Lemma 7.10.} \textsl{The constant }$\kappa $\textsl{\ in
the formula (6.9) is }$2(-1)^{k(n)-1}$\textsl{.}

\bigskip

\noindent \textbf{Proof.} Let $D$ be the ideal on $H^{\ast }(B_{U(k)})=%
\mathbb{Z}[c_{1},\cdots ,c_{k}]$ generated by the $c_{r}$ with $r\geq 2$,
and consider the ring map $e:H^{\ast }(B_{U(k)})\rightarrow \mathbb{Z}$
defined by

\begin{quote}
$e(c_{1})=1$, $e(c_{r})=0$, $2\leq r\leq k$.
\end{quote}

\noindent That is, for every $u\in H^{2k}(B_{U(k)})$,

\begin{quote}
i) $u=e(u)\cdot c_{1}^{k}+h(u)$ with $h(u)\in D$;

ii) $e(u)=0$ if and only if $u\in D$.
\end{quote}

\noindent In particular, applying $e$ to the relations

\begin{quote}
$2\alpha _{r+1}+\alpha _{r}^{2}=f_{r}(g_{1},\cdots ,g_{k-1})$ (by Theorem D'
and (7.4))
\end{quote}

\noindent on $H^{\ast }(B_{U(k)})$ we get from $e(\alpha _{1})=2$ (by
Example 7.9) that

\begin{enumerate}
\item[(7.13)] $e(\alpha _{r})=2(-1)^{r-1}$.
\end{enumerate}

On the other hand, since $B_{\lambda }^{\ast }(\overline{\theta }_{n})$ is
the Euler class of the induced bundle $B_{\lambda _{n}}^{\ast }\eta _{n}$ on 
$B_{U(k)}$ we have

\begin{quote}
$e(B_{\lambda }^{\ast }(\overline{\theta }_{n}))=1$ (i.e. $B_{\lambda
}^{\ast }(\overline{\theta }_{n})=c_{1}^{2^{k(n)}}+\alpha $ for some $\alpha
\in D$).
\end{quote}

\noindent Thus, applying the ring map $e\circ B_{\lambda }^{\ast }$ to the
equation (6.8) and noting that $B_{\lambda }^{\ast }(\varphi (\overline{b}%
_{n}))\in D$ (since $\overline{b}_{n}\in A_{1}$), we obtain $%
2(-1)^{k(n)-1}=\kappa $ as required.$\square $

\bigskip

\noindent \textbf{Remark 7.10. }For a $CW$ complex $X$ the maps $B_{\lambda
_{0}}$, $B_{\lambda }$ and $B_{\lambda ^{c}}$ induce, respectively, the
correspondences between homotopy sets

\begin{quote}
$B_{\lambda _{0\ast }}:[X,B_{U(k)}]\rightarrow \lbrack X,B_{SO(2k)}]$ by $%
B_{\lambda _{0}\ast }[g]=[B_{\lambda _{0}}\circ g]$,

$B_{\lambda _{\ast }}:[X,B_{U(k)}]\rightarrow \lbrack X,B_{Spin(2k)}]$ by $%
B_{\lambda _{\ast }}[g]=[B_{\lambda }\circ g]$,

$B_{\lambda _{\ast }^{c}}:[X,B_{U^{c}(k)}]\rightarrow \lbrack
X,B_{Spin^{c}(2k)}]$ by $B_{\lambda _{\ast }^{c}}[g]=[B_{\lambda ^{c}}\circ
g]$,
\end{quote}

\noindent in which $B_{\lambda _{0\ast }}$ is well known to be \textsl{the
real reduction} on the complex bundles \cite[p.155]{MS}. Likewise, the maps $%
B_{\lambda _{\ast }}$ and $B_{\lambda _{\ast }^{c}}$ can be regarded as the 
\textsl{spin} and the \textsl{spin}$^{c}$\textsl{\ reduction} of complex
bundles, respectively. In this connection, the formulae in Theorem 7.8
express the Spin characteristic classes (see in \S 9) of the spin reduction
of a complex vector bundle by its Chern classes.

The map $B_{\lambda }$ fits also into the fibration

\begin{quote}
$Spin(2k)/U(k)\hookrightarrow B_{U(k)}\rightarrow B_{Spin(2k)}$
\end{quote}

\noindent that is also of geometric significances \cite{D,D1}: the fiber
manifold $Spin(2k)/U(k)$ acts as the Grassmanian of complex structures on
the $2k$ dimensional Euclidean space $\mathbb{R}^{2k}$, and can be
identified with the classifying space of the complex $k$ bundles whose real
(resp. $Spin$) reductions are trivial.$\square $

\section{The Weyl invariants of the groups $Spin(n)$}

Let $G$ be a compact connected Lie group with a maximal torus $T$, and the
Weyl group $W=N_{G}(T)/T$. The canonical $W$ action on $T$ induces an action
on the integral cohomology $H^{\ast }(B_{T})$. Denote by $H^{\ast
}(B_{T})^{W}$ the subring consisting of all the $W$ invariants, and let $%
B_{t}:B_{T}\rightarrow B_{G}$ be induced by the inclusion $t:T\rightarrow G$%
. A classical result of Borel \cite{B2} states that

\bigskip

\noindent \textbf{Lemma 8.1.} \textsl{The ring map }$B_{t}^{\ast }:H^{\ast
}(B_{G})\rightarrow H^{\ast }(B_{T})$\textsl{\ annihilates the torsion ideal 
}$\tau (B_{G})$\textsl{, and induces an injection }$H^{\ast }(B_{G})/\tau
(B_{G})\rightarrow H^{\ast }(B_{T})^{W}$\textsl{.}$\square $

\bigskip

The fundamental problem of the invariant theory of Weyl groups is to present
the ring $H^{\ast }(B_{T})^{W}$ by explicit generators and relations \cite%
{Gu,We}. A closely related problem in topology is to decide the subring $%
\func{Im}B_{t}^{\ast }\subseteq H^{\ast }(B_{T})^{W}$. For the Weyl group $W$
of a semi--simple Lie group $G$ Chevalley has shown that

\begin{quote}
$H^{\ast }(B_{T})^{W}\otimes \mathbb{Z}_{0}=H^{\ast }(B_{G})\otimes \mathbb{Z%
}_{0}=\mathbb{Z}_{0}[P_{1},\cdots ,P_{n}]$, $n=\dim T$,
\end{quote}

\noindent where \textsl{the basic (rational) polynomial} \textsl{invariants }%
$\left\{ P_{1},\cdots ,P_{n}\in H^{\ast }(B_{T})^{W}\otimes \mathbb{Z}%
_{0}\right\} $ have been made explicitly by Mehta \cite{Me}. In addition,
for a prime $p$ methods to calculate the algebras $H^{\ast
}(B_{T})^{W}\otimes \mathbb{Z}_{p}$ have been developed from the
perspectives of algebraic topology \cite{KM,Sm}, combinatorics \cite{St},
and algorithm \cite{Bt}. However, apart from Borel's classical result \cite[%
\S 3. Examples]{Fe1}

\begin{quote}
$\func{Im}B_{t}^{\ast }=H^{\ast }(B_{T})^{W}=H^{\ast }(B_{G})$ for $G=SU(n)$
or $Sp(n)$,
\end{quote}

\noindent complete information on the ring $H^{\ast }(B_{T})^{W}$ are not
known for most other semi--simple Lie groups. In particular, the types of
the Weyl groups of $G=Spin(n)$ consist of $B_{k}$ ($n=2k+1$) or $D_{k}$ ($%
n=2k$), and partial information has been obtained by Borel, Feshbach,
Totaro, Benson and Wood \cite{B2,BW,Fe1,T}.

For $G=Spin(n)$ our approach to $H^{\ast }(B_{T})^{W}$ begins with the
subring $\func{Im}B_{t}^{\ast }$. For an integer $n>6$ (see Remark 1.1) we
set $k=\left[ \frac{n}{2}\right] $, and let $h$ be the inclusion of the
diagonal subgroup $T=U(1)\times \cdots \times U(1)$ ($k$ copies) into $U(k)$%
. Then a convenient maximal torus on $Spin(n)$ is

\begin{quote}
$t=\lambda \circ h:T\rightarrow U(k)\rightarrow Spin(2k)\subseteq Spin(n)$,
\end{quote}

\noindent where $\lambda $ is the inclusion given in table (7.1).
Furthermore, with respect to the canonical presentation

\begin{enumerate}
\item[(8.1)] $H^{\ast }(B_{T})=\mathbb{Z}[x_{1},\cdots ,x_{k}]$, $\deg
x_{i}=2$,
\end{enumerate}

\noindent the ring map $B_{h}^{\ast }:$ $H^{\ast }(B_{U(k)})\rightarrow
H^{\ast }(B_{T})$ is given by

\begin{quote}
$B_{h}^{\ast }(c_{r})=e_{r}(x_{1},\cdots ,x_{k})$, $1\leq r\leq k$,
\end{quote}

\noindent where $e_{r}$ is the $r^{th}$ elementary symmetric function in the 
$x_{i}$'s. It follows that $B_{h}^{\ast }$ carries $H^{\ast }(B_{U(k)})$\
isomorphically onto the subring $Sym[x_{1},\cdots ,x_{k}]$\ of symmetric
functions, while $B_{\lambda }^{\ast }$\ induces an injection from the
quotient ring $H^{\ast }(B_{Spin(n)})/\tau (B_{Spin(n)})$ into $H^{\ast
}(B_{U(k)})$. Thus, applying the ring map $B_{t}^{\ast }$ to the formula
(6.10) of the ring $H^{\ast }(B_{Spin(n)})$ we obtain by Theorem 7.8 the
following characterization of the subring $\func{Im}B_{t}^{\ast }\subseteq
H^{\ast }(B_{T})^{W}$.

\bigskip

\noindent \textbf{Theorem 8.2.} \textsl{The subring }$\func{Im}B_{t}^{\ast }$%
\textsl{\ has the presentations:}

\begin{enumerate}
\item[(8.2)] $\func{Im}B_{t}^{\ast }=\left\{ 
\begin{tabular}{l}
$\mathbb{Z}[g_{1},\cdots ,g_{k-1},c_{k},\alpha _{1},\cdots ,\alpha
_{k(n)-1},B_{t}^{\ast }(\overline{\theta }_{n})]/D_{n}$ \textsl{if} $n=2k$;
\\ 
$\mathbb{Z}[g_{1},\cdots ,g_{k},\alpha _{1},\cdots ,\alpha
_{k(n)-1},B_{t}^{\ast }(\overline{\theta }_{n})]/D_{n}$ \textsl{if} $n=2k+1$,%
\end{tabular}%
\right. $
\end{enumerate}

\noindent \textsl{where }$D_{n}$\textsl{\ denotes the ideal generated by the
following relations}

\begin{quote}
\textsl{i)} $2\alpha _{1}-g_{1}$\textsl{,} $2\alpha _{r+1}+\alpha
_{r}^{2}-f_{r}(g_{1},\cdots ,g_{k})$\textsl{,} $1\leq r\leq k(n)-2$\textsl{,}

\textsl{ii)} $(-1)^{k(n)-1}\cdot 4\cdot B_{t}^{\ast }(\overline{\theta }%
_{n})+\alpha _{k(n)-1}^{2}-\beta _{n}$\textsl{,}
\end{quote}

\noindent \textsl{where }$\beta _{n}:=B_{t}^{\ast }(\overline{b}_{n})$%
\textsl{,} \textsl{and where the invariants }$g_{r}$\textsl{,} $\alpha _{r}$ 
\textsl{and} $f_{r}(g_{1},\cdots ,g_{k})$ \textsl{are given by formulae
(7.10), (7.11) and (7.12), respectively.}$\square $

\bigskip

In term of (8.2) define $\func{Im}\overline{B}_{t}^{\ast }\subset $\ $\func{%
Im}B_{t}^{\ast }$ to be the subring generated by

\begin{quote}
$g_{2},\cdots ,g_{\left[ \frac{n-1}{2}\right] },\alpha _{1},\cdots ,$ $%
\alpha _{k(n)-2}$,\ and $c_{k}$\ if $n=2k$.
\end{quote}

\noindent Then, for the degree reason, the relation $2\alpha
_{k(n)-1}+\alpha _{k(n)-2}^{2}=f_{k(n)-2}(g_{1},\cdots ,g_{k})$ on $\func{Im}%
B_{t}^{\ast }$ implies that

\bigskip

\noindent \textbf{Corollary 8.3. }\textsl{The quotient group }$\func{Im}%
B_{t}^{2^{k(n)}}/\func{Im}\overline{B}_{t}^{2^{k(n)}}$\textsl{\ is
isomorphic to }$\mathbb{Z}_{2}$\textsl{\ with generator }$\alpha _{k(n)-1}$%
\textsl{.}$\square $

\bigskip

For $G=Spin(n)$ the extension problem from $\func{Im}B_{t}^{\ast }$ to the
ring $H^{\ast }(B_{T})^{W}$ was raised by Borel \cite[1954]{B2}, studied by
Feshbach \cite[1981]{Fe1}, and has been solved by Benson and Wood in the
remarkable work \cite[1995]{BW}. Efforts to bring the relevant constructions
and calculations taking place in \cite{BW} into our context gives rise to
the following results.

\bigskip

\noindent \textbf{Theorem 8.4 (Benson and Wood).} \textsl{Assume that }$%
G=Spin(n)$ \textsl{with} $n>6$\textsl{.}

\textsl{i) If }$n\neq 3,4,5\func{mod}8$\textsl{\ then }$\func{Im}B_{t}^{\ast
}=H^{\ast }(B_{T})^{W}$\textsl{.}

\textsl{ii) If} $n$\textsl{\ }$=3,4,5\func{mod}8$\textsl{,\ then }$H^{\ast
}(B_{T})^{W}$ \textsl{is} \textsl{generated by its subring }$\func{Im}%
B_{t}^{\ast }$\textsl{, together with an additional invariant }$\omega
_{k(n)-1}\in H^{2^{k(n)}}(B_{T})^{W}$\textsl{\ that is related to the known
invariants }$B_{t}^{\ast }(\overline{\theta }_{n})$ \textsl{and} $\alpha
_{k(n)-1}$ \textsl{by the relations}

\begin{quote}
\textsl{a) }$\omega _{k(n)-1}^{2}=B_{\lambda }^{\ast }(\overline{\theta }%
_{n})$\textsl{; }

\textsl{b) }$2\omega _{k(n)-1}-\alpha _{k(n)-1}=l_{n}$\textsl{, }$l_{n}\in 
\func{Im}\overline{B}_{t}^{2^{k(n)}}$\textsl{.}
\end{quote}

\noindent \textbf{Proof. }With respect to $n=2k+1$ or $n=2k$ consider the
sequence of cohomology classes introduced in \cite[\S 4]{BW}

\begin{quote}
$\left\{ \eta _{j}\in H^{\ast }(B_{T})\text{, }1\leq j\leq k\right\} $ or $%
\left\{ \mu _{j}\in H^{\ast }(B_{T})\text{, }1\leq j\leq k\right\} $.
\end{quote}

\noindent Then by \cite[Table 2]{BW} we have

\begin{enumerate}
\item[(8.3)] $B_{\lambda }^{\ast }(\overline{\theta }_{n})=\eta _{k(n)}$ or $%
\mu _{k(n)\text{ }}$ if $n\equiv 3,5\func{mod}8$ or $n\equiv 4\func{mod}8$.
\end{enumerate}

\noindent Let us put

\begin{quote}
$\omega _{k(n)-1}:=\eta _{k(n)-1}$ if $n\equiv 3,5\func{mod}8$, or $\mu
_{k(n)-1}$ if $n\equiv 4\func{mod}8$.
\end{quote}

\noindent Then Benson and Wood \cite[Proposition 4.1]{BW} has shown that

\begin{enumerate}
\item[(8.4)] $\omega _{k(n)-1}\in H^{\ast }(B_{T})^{W}$ and $\omega
_{k(n)-1}^{2}=B_{t}^{\ast }(\overline{\theta }_{n})$.
\end{enumerate}

\noindent Now, apart from for the relation b) all the statements of the
theorem are verified by comparing \cite[Theorem 7.1]{BW} with \cite[Theorem
10.2]{BW}.

Consider the sequence $\left\{ q_{r},r\geq 1\right\} $ on $H^{\ast
}(B_{T})^{W}$ constructed in the proof of \cite[Proposition 3.3]{BW}. By 
\cite[Corollary 7.2]{BW}

\begin{enumerate}
\item[(8.5)] $2\omega _{k(n)-1}-q_{k(n)-1}=l_{n}$ for some $l_{n}\in \func{Im%
}\overline{B}_{t}^{2^{k(n)}}$.
\end{enumerate}

\noindent On the other hand, by ii) and iii) of \cite[Theorem 10.2]{BW} the
class $q_{k(n)-1}$ generates also the quotient group

\begin{quote}
$\func{Im}B_{t}^{2^{k(n)}}/\func{Im}\overline{B}_{t}^{2^{k(n)}}=\mathbb{Z}%
_{2}$.
\end{quote}

\noindent Comparing this with Corollary 8.3 we can replace in (8.5) the
class $q_{k(n)-1}$ by our $\alpha _{k(n)-1}$ to obtain the desired relation
b).$\square $

\bigskip

Assume that $n$\textsl{\ }$=3,4,5\func{mod}8$. The relations a) and b) of
Theorem 8.4 imply that the generators $B_{t}^{\ast }(\overline{\theta }_{n})$%
\textsl{\ }and\textsl{\ }$\alpha _{k(n)-1}$ of $\func{Im}B_{t}^{\ast }$ can
be expressed as polynomials in the elements

\begin{quote}
$\omega _{k(n)-1}$, $g_{2},\cdots ,g_{\left[ \frac{n-1}{2}\right] },\alpha
_{1},\cdots ,\alpha _{k(n)-2}$, and $c_{k}$\ if $n=2k$.
\end{quote}

\noindent In addition, combining the relation b) with the relation on $\func{%
Im}B_{t}^{\ast }$

\begin{quote}
$2\alpha _{k(n)-1}+\alpha _{k(n)-2}^{2}=f_{k(n)-2}(g_{1},\cdots ,g_{k})$ (by
Theorem 8.2)
\end{quote}

\noindent one gets

\begin{quote}
$4\omega _{k(n)-1}-\alpha _{k(n)-2}^{2}=\varepsilon _{n}$ for some $%
\varepsilon _{n}\in \func{Im}\overline{B}_{t}^{\ast }$.
\end{quote}

\noindent Thus, putting Theorems 8.2 and Theorem 8.4 together we obtain that

\bigskip

\noindent \textbf{Theorem 8.5.} \textsl{Assume that }$G=Spin(n)$ \textsl{with%
} $n>6$\textsl{\ (see Remark 1.1).} \textsl{Then }

\begin{enumerate}
\item[(8.6)] $H^{\ast }(B_{T})^{W}=\left\{ 
\begin{tabular}{l}
$\func{Im}B_{t}^{\ast }\text{ \textsl{if} }n\QTR{sl}{\ }\neq 3,4,5\func{mod}8%
\text{;}$ \\ 
$\func{Im}\overline{B}_{t}^{\ast }\otimes \mathbb{Z}[\omega
_{k(n)-1}]/\left\langle h_{n}\right\rangle \text{ \textsl{if} }n\QTR{sl}{\ }%
\equiv 3,4,5\func{mod}8\text{,}$%
\end{tabular}%
\right. $
\end{enumerate}

\noindent \textsl{where }$h_{n}=4\cdot \omega _{k(n)-1}-\alpha
_{k(n)-2}^{2}-\varepsilon _{n}$\textsl{\ with }$\varepsilon _{n}\in \func{Im}%
\overline{B}_{t}^{\ast }.\square $

\bigskip

\noindent \textbf{Remarks 8.6. }Theorems 8.2 and 8.5 present both of the
rings $\func{Im}B_{t}^{\ast }$ and $H^{\ast }(B_{T})^{W}$ by the explicit
generators (e.g. (7.10) and (7.11))

\begin{quote}
$g_{1},\cdots ,g_{k-1},c_{k},\alpha _{1},\cdots ,\alpha
_{k(n)-1},B_{t}^{\ast }(\overline{\theta }_{n})$
\end{quote}

\noindent together with the invariant $\omega _{k(n)-1}$ given by Benson and
Wood \cite[\S 4]{BW}.

In \cite[Theorem 7.1]{BW} Benson and Wood obtained a presentation of the
ring $H^{\ast }(B_{T})^{W}$ without specifying the relations among their
generators. In our approach the recurrence relations

\begin{quote}
$2\alpha _{1}-g_{1}$, $2\alpha _{r+1}+\alpha _{r}^{2}-f_{r}(g_{1},\cdots
,g_{k}),$ $1\leq r\leq k(n)-2$,
\end{quote}

\noindent are originated from property ii) of Theorem C', which are useful
to produce the sequence $\{\alpha _{1},\alpha _{2},\cdots \}$ of invariants
from the initial one $\alpha _{1}=-c_{2}+2c_{1}^{2}$, see Examples 7.6 and
7.9.$\square $

\section{The Spin characteristic classes}

For the groups $SO=\cup _{n=2}^{\infty }SO(n)$ and $Spin=\cup _{n=2}^{\infty
}Spin(n)$ in the stable range we have by formulae (2.3) and (6.10) that

\begin{quote}
$H^{\ast }(B_{SO})=\mathbb{Z}[p_{1},p_{2},\cdots ]\oplus \tau (B_{SO})$ with 
$2\cdot \tau (B_{SO})=0$,

$H^{\ast }(B_{Spin})=\overline{\pi }^{\ast }H^{\ast }(B_{SO})\otimes \mathbb{%
Z}[\overline{q}_{1},\overline{q}_{2},\cdots ]/K_{\infty }$,
\end{quote}

\noindent respectively, where the Euler classes $\overline{\pi }^{\ast
}(e_{n})$, $\overline{\theta }_{n}$ are disappeared at $n=\infty $. In view
of these formulae we introduce the sequence $\left\{ Q_{k},\text{ }k\geq
1\right\} $, $\deg Q_{k}=4k$, of integral cohomology classes on $B_{Spin}$
by setting

\begin{quote}
$Q_{k}:=\left\{ 
\begin{tabular}{l}
$\overline{\pi }^{\ast }p_{k}$ if $k>1$ is not a power of $2$; \\ 
$\overline{q}_{r}$ if $k=2^{r}$, $r\geq 0$.%
\end{tabular}%
\right. $
\end{quote}

\noindent Then the relation ii) of Theorem C' implies the formulae

\begin{quote}
$2Q_{1}=\overline{\pi }^{\ast }p_{1}$, $2Q_{2^{r}}+Q_{2^{r-1}}^{2}=\overline{%
\pi }^{\ast }f(w_{2}^{(r)})$
\end{quote}

\noindent in which by Example 2.4 and by the formula (2.4) of $f$

\begin{quote}
$f(w_{2}^{(r)})=p_{2^{r}}+p_{1}p_{2^{r}-2}+\cdots
+p_{2^{r-1}-2}p_{2^{r-1}+2}+$ higher terms.
\end{quote}

\noindent This implies that\ every $2$--power Pontryagin class $\overline{%
\pi }^{\ast }p_{2^{r}}$\ can be expressed as a polynomial in the $Q_{k}$'s,
plus some torsion element. Summarizing, taking $n=\infty $ in Theorem D' we
can eliminate

\begin{quote}
$\overline{\pi }^{\ast }(e_{n})$, $\overline{\theta }_{n}$, $\overline{\pi }%
^{\ast }p_{2^{r}}$ and $2Q_{2^{r}}+Q_{2^{r-1}}^{2}=\overline{\pi }^{\ast
}f(w_{2}^{(r)})$
\end{quote}

\noindent form the sets of generators and relations to obtain the following
result from Theorem C', as well as the formula (6.6) of $\tau (B_{Spin(n)})$.

\bigskip

\noindent \textbf{Theorem 9.1.} \textsl{The integral cohomology of }$%
B_{Spin} $\textsl{\ has the presentation}

\begin{enumerate}
\item[(9.1)] $H^{\ast }(B_{Spin})=\mathbb{Z}[Q_{1},Q_{2},Q_{3},\cdots
]\oplus \overline{\pi }^{\ast }\tau (B_{SO})$\textsl{, }$\deg Q_{k}=4k$%
\textsl{,}
\end{enumerate}

\noindent \textsl{where} \textsl{the generators }$Q_{k}$\textsl{\ are
characterized uniquely by the following properties:}

\textsl{i) if }$k>1$\textsl{\ is not a power of }$2$\textsl{, then }$Q_{k}=%
\overline{\pi }^{\ast }p_{k}$\textsl{;}

\textsl{ii) if }$k=2^{r}$\textsl{\ with }$r\geq 0$ \textsl{then}

\begin{enumerate}
\item[(9.2)] $\rho _{2}(Q_{k})=\overline{\pi }^{\ast }(w_{2}^{(k+1)})$%
\textsl{,}

\item[(9.3)] $2Q_{1}=\overline{\pi }^{\ast }(p_{1})$\textsl{,} $%
2Q_{2k}+Q_{k}^{2}=\overline{\pi }^{\ast }f(w_{2}^{(k+1)})$\textsl{, }$k\geq
1 $\textsl{,}
\end{enumerate}

\noindent \textsl{In particular, the cup product between the free part and
the torsion ideal of the ring }$H^{\ast }(B_{Spin})$ \textsl{is given by}

\begin{center}
$Q_{k}\cup \overline{\pi }^{\ast }(\delta _{2}(x))=\left\{ 
\begin{tabular}{l}
$\overline{\pi }^{\ast }(\delta _{2}(x\cup w_{2k}^{2}))\text{ }$\textsl{if }$%
k>1$\textsl{\ is not a power of} $2$\textsl{;} \\ 
$\overline{\pi }^{\ast }(\delta _{2}(x\cup w_{2}^{(r+1)}))\text{ \textsl{if} 
}k=\text{ }2^{r}\text{.}\square $%
\end{tabular}%
\right. $
\end{center}

In the formula (9.1) the torsion ideal $\overline{\pi }^{\ast }\tau (B_{SO})$
is fashioned from $\tau (B_{SO})$, hence contributes nothing essentially
new. It is the generators $\{Q_{k},k\geq 1\}$ of the free part, together
with their uniqueness property, are just what requested for us to define the
characteristic classes for the Spin vector bundles.

Let $\xi $ be an oriented real bundle over a connected $CW$--complex $X$,
induced by a map $f_{\xi }$ from $X$ to $B_{SO}$, and suppose that $%
w_{2}(\xi )=0$ (i.e. $\xi $ is Spin). Then $f_{\xi }$ can be factored into a
composition

\begin{quote}
$X\overset{g_{\xi }}{\rightarrow }B_{Spin}\overset{\overline{\pi }}{%
\rightarrow }B_{SO}$,
\end{quote}

\noindent where the map $g_{\xi }$ is unique up to homotopy.

\bigskip

\noindent \textbf{Definition 9.2.} The cohomology class $q_{k}(\xi ):=g_{\xi
}^{\ast }(Q_{k})\in H^{4k}(X)$, $k\geq 1$, (resp. the sum $q(\xi
):=1+q_{1}(\xi )+q_{2}(\xi )+\cdots \in H^{\ast }(X)$) is called the $k^{th}$
\textsl{Spin\ characteristic class} (resp. \textsl{the total Spin\
characteristic class}) of the bundle $\xi $.$\square $

\bigskip

The Spin characteristic classes so defined possesses the naturality property.

\bigskip

\noindent \textbf{Corollary 9.3.} \textsl{For any map }$g:Y\rightarrow X$%
\textsl{\ between CW--complexes one has}

\begin{quote}
$q_{k}(g^{\ast }\xi )=g^{\ast }q_{k}(\xi )$\textsl{, }$k\geq 1$\textsl{.}
\end{quote}

\noindent \textsl{In particular, the generator }$Q_{2^{r}}$\textsl{\ of }$%
H^{\ast }(B_{Spin})$ \textsl{restricts to} \textsl{the generator }$\overline{%
q}_{r}$\textsl{\ of }$H^{\ast }(B_{Spin(n)})$\textsl{\ (see (6.10)) via the
inclusion }$B_{Spin(n)}\subset B_{Spin}$\textsl{.}$\square $

\bigskip

\noindent \textbf{Example 9.4. }The relations (9.2) and (9.3) have
non--trivial implications. Let $\xi $ be a spin vector bundle over a space $%
X $, $\dim \xi =n$, and let

\begin{quote}
$w(\xi )=1+w_{4}(\xi )+$ $\cdots +w_{n}(\xi )\in H^{\ast }(X;\mathbb{Z}_{2})$
or

$p(\xi )=1+p_{1}(\xi )+p_{2}(\xi )+$ $\cdots +p_{\left[ \frac{n-1}{2}\right]
}(\xi )\in H^{\ast }(X)$
\end{quote}

\noindent be its total Stiefel--Whitney or Pontryagin classes, respectively.
Then (9.2) turns out

\begin{quote}
$\rho _{2}(q_{1}(\xi ))=w_{4}(\xi )$,

$\rho _{2}(q_{2}(\xi ))=w_{8}(\xi ),$

$\rho _{2}(q_{4}(\xi ))=w_{16}(\xi )+w_{4}(\xi )w_{12}(\xi )+w_{6}(\xi
)w_{10}(\xi )+w_{4}(\xi )w_{6}^{2}(\xi )$, $\cdots $
\end{quote}

\noindent by Example 2.4. Since these polynomials admit integral lifts,
applying $Sq^{1}$ to both sides yields the following universal relations
among the Stiefel--Whitney classes of a spin bundle $\xi $:

\begin{quote}
$w_{5}(\xi )=0$, $w_{9}(\xi )=0$,

$w_{17}(\xi )+w_{4}(\xi )w_{13}(\xi )+w_{7}(\xi )w_{10}(\xi )+w_{6}(\xi
)w_{11}(\xi )=0$, $\cdots $.
\end{quote}

Ignoring the elements of order $2$ the relation (9.3) allows us to express
the Pontryagin classes $p_{i}(\xi )$'s as polynomials in the Spin ones $%
q_{i}(\xi )$'s, such as

\begin{enumerate}
\item[(9.4)] $p_{1}=2q_{1}$, $p_{2}=2q_{2}+q_{1}^{2}$, $p_{3}=q_{3}$, $%
p_{4}=2q_{4}+q_{2}^{2}-2q_{1}q_{3}$, $\cdots $,
\end{enumerate}

\noindent by Example 2.4 and formula (2.4). In Section \S 10 these
transition functions will be applied to simplify various formulae of spin
manifolds.$\square $

\bigskip

\noindent \textbf{Remark 9.5.} The spinors were first discovered by E.
Cartan in 1913 in his investigations of the representation theory of
topological groups, and had subsequently found significant and wide
applications to geometry and mathematical physics. However, a precise
definition of spin structure was possible only after the notion of fiber
bundle was introduced. Notably, Haefliger (1956) proved\ that the second
Stiefel--Whitney class $w_{2}(M)$ is the only obstruction to the existence
of a spin structure on an oriented Riemannian manifold $M$. This was soon
extended by Borel and Hirzebruch (1958) to the cases of vector bundles over
CW--complexes.

The idea of \textsl{Spin characteristics} was initiated by Thomas. He \cite[%
Theorem (1.2)]{Th} described the integral cohomology $H^{\ast }(B_{Spin})$
using a squence $\{Q_{j}\}$ of generators, but that is subject to two
sequences $\{\Phi _{j}\}$ and $\{\Psi _{j}\}$ of indeterminacies, where he
asked also for axioms with geometric significance by which the uniqueness of
such a sequence $\{Q_{j}\}$ can be secured.

Granted with the two sequences $\{w_{2},w_{2}^{(1)},\cdots \}$ and $%
\{f(w_{2}),f(w_{2}^{(1)}),\cdots \}$ of cohomology classes depending on the
only obstruction $w_{2}$ to the existence of spin structure, Theorem C
(resp. Theorem C') amounts to an axiomatic characterization of our
generators $\left\{ q_{r},r\geq 0\right\} $ on $H^{\ast }(B_{Spin^{c}(n)})$
(resp. $\left\{ \overline{q}_{r},r\geq 0\right\} $ on $H^{\ast
}(B_{Spin(n)}) $). Not surprisingly, the Spin characteristic classes so
obtained are better adapted with topics of spin geometry, see in Section \S %
10.$\square $

\section{Applications to spin geometry}

For a smooth manifold $M$ its total Pontryagin class (resp. total Stiefel
Whitney class) is defined to be that of the tangent bundle $TM$ of $M$, and
is denoted by

\begin{quote}
$p(M):=1+p_{1}+\cdots +p_{k}$, $k=\left[ \frac{n}{4}\right] $

(resp. $w(M):=1+w_{1}+\cdots +w_{n}$, $n=\dim M$).
\end{quote}

\noindent Similarly, if $M$\ is spin (i.e. $w_{2}(M)=0$), then its total
Spin characteristic class is also defined, and will be written as

\begin{quote}
$q(M):=1+q_{1}+\cdots +q_{k}$, $k=\left[ \frac{\dim M}{4}\right] $.
\end{quote}

In the topological approach to spin geometry, the Spin characteristic
classes can play roles that may not be replaced by the regular
characteristic classes. In this section we provide such initial evidences.
Subject to the main theme of this paper the examples and calculations will
be restricted to relatively lower dimensional cases.

\bigskip

\textbf{10.1. The tangent invariants of Spin manifolds.} According to C.T.C.
Wall \cite[Theorem 5]{Wa1} for each integer $b$ there exists a unique $6$%
--dimensional smooth manifold $\mathbb{C}P_{b}^{3}$ homotopy equivalent to
the $3$--dimensional complex projective space $\mathbb{C}P^{3}$, whose first
Pontryagin class is $p_{1}=4(1+6b)x^{2}$, where $x$ denotes a generator of $%
H^{2}(\mathbb{C}P_{b}^{3})=\mathbb{Z}$. Since the total Stiefel--Whitney
class of a manifold is a homotopy invariant we must have $w(\mathbb{C}%
P_{b}^{3})=w(\mathbb{C}P^{3})=1$. In particular, the manifold $\mathbb{C}%
P_{b}^{3}$ is spin with

\begin{quote}
$q_{1}(\mathbb{C}P_{b}^{3})=\frac{1}{2}p_{1}(\mathbb{C}%
P_{b}^{3})=2(1+6b)x^{2}$ (by (9.4)).
\end{quote}

Let $\pi :M_{b}^{7}\rightarrow \mathbb{C}P_{b}^{3}$ be the oriented circle
bundle on $\mathbb{C}P_{b}^{3}$ with Euler class $e=4(1+6b)x$. From the
Gysin sequence of $\pi $ \cite[p.157]{MS} one sees that the groups $%
H^{2r}(M_{b}^{7})$ is cyclic of order $4(1+6b)$ with generators $\pi ^{\ast
}(x^{r})$, where $r=1,2,3$. In view of the decomposition $TM_{b}^{7}=\pi
^{\ast }T\mathbb{C}P_{b}^{3}\oplus \varepsilon $ on the tangent bundle of $%
M_{b}^{7}$ we get by the naturality of the characteristic classes that

\begin{quote}
$w(M_{b}^{7})=1$,\textsl{\ }$p(M_{b}^{7})=1$,\textsl{\ }but $%
q_{1}(M_{b}^{7})=2(1+6b)\pi ^{\ast }(x^{2})\neq 0$,
\end{quote}

\noindent where $\varepsilon $ denotes the $1$--dimensional trivial bundle
on $M_{b}^{7}$. This shows that:

\bigskip

\noindent \textbf{Proposition 10.1.} \textsl{The family }$\{M_{b}^{7},b\in 
\mathbb{Z}\}$ \textsl{of} $7$\textsl{--dimensional} \textsl{smooth} \textsl{%
Spin manifolds satisfies that }$w(M_{b}^{7})=1$\textsl{, }$p(M_{b}^{7})=1,$%
\textsl{\ but }$q(M_{b}^{7})\neq 1$\textsl{.}$\square $

\bigskip

\textbf{10.2. The integral lifts of Wu--classes.} For a sufficiently large $%
n $ let $v_{r}\in H^{r}(B_{Spin(n)};\mathbb{Z}_{2})$ be the $r^{th}$
Wu--class of the canonical real $n$--bundle on $B_{Spin(n)}$. It is well
known that $v_{r}=0$ unless $r\equiv 0\func{mod}4$, and that all the classes 
$v_{4k}$ admit integral lifts (see \cite{ABP} or \cite[Lemma E1]{HS}). In 
\cite{HS} Hopkins and Singer constructed a stable exponential characteristic
class $v_{t}^{Spin}$ with values in the integral cohomology $H^{\ast
}(B_{Spin(n)})$, whose mod $2$--reduction is the total Wu--class.
Rationally, in terms of Pontryagin classes, the first four terms are

\begin{quote}
$v_{4}^{Spin}=-\frac{1}{2}p_{1}$,

$v_{8}^{Spin}=\frac{1}{2^{3}}(20p_{2}-9p_{1}^{2})$,

$v_{12}^{Spin}=-\frac{1}{2^{4}}(80p_{3}+60p_{1}p_{2}-17p_{1}^{3})$,

$v_{16}^{Spin}=\frac{1}{2^{7}}(2^{6}\cdot 29p_{4}-2^{4}\cdot
33p_{2}^{2}+2^{3}\cdot 147p_{1}^{2}p_{2}-277p_{1}^{4})$.
\end{quote}

\noindent Substituting the Pontryagin classes by the Spin characteristic
classes using the transition (9.4) shows that

\bigskip

\noindent \textbf{Proposition 10.2.} \textsl{In term of the Spin
characteristic classes, the first four Wu--classes }$v_{4k}$, $i=1,2,3$ 
\textsl{or} $4,$\textsl{\ admit the integral lifts:}

\begin{quote}
$\widetilde{v}_{4}=q_{1}$\textsl{, }

$\widetilde{v}_{8}=q_{2}$\textsl{, }

$\widetilde{v}_{12}=q_{3}+q_{1}q_{2}+q_{1}^{3}$\textsl{,}

$\widetilde{v}_{16}=q_{4}+q_{1}q_{3}+q_{1}^{2}q_{2}$\textsl{.}$\square $
\end{quote}

For possible applications of these formulae in geometry we refer to Wilson 
\cite{Wi}, Landweber and Stong \cite{LS}, where the following subtle
relations are found for the spin manifolds of dimension $8k+2$

\begin{enumerate}
\item[(10.1)] $Sq^{3}v_{4k}=0$, $w_{4}w_{8k-2}=v_{4k}Sq^{2}v_{4k}$.
\end{enumerate}

\textbf{10.3. The Eells--Kuiper invariant.} For a closed, smooth, oriented $%
(4k-1)$ manifold $M$ furnished with a spin coboundary $W$ that satisfies so
called $\mu $--conditions, Eells and Kuiper \cite{EK} introduced a
differential invariant $\mu _{k}(M)$ of $M$ in terms of Pontryagin numbers
and the signature $\sigma $ of the coboundary $W$. For the small values $%
k=2,3,4$ the formulae of $\mu _{k}(M)$ reads

\begin{quote}
$\mu _{2}(M)\equiv \frac{(p_{1}^{2}-4\sigma )[W]}{2^{7}\cdot 7}\func{mod}1$

$\mu _{3}(M)\equiv \frac{(4p_{1}p_{2}-3p_{1}^{3}-24\sigma )[W]}{2^{11}\cdot
3\cdot 31}\func{mod}1$

$\mu _{4}(M)\equiv \frac{%
(12096p_{1}p_{3}+5040p_{2}^{2}-22680p_{1}^{2}p_{2}+9639p_{1}^{4}-181440%
\sigma )[W]}{2^{15}\cdot 3^{4}\cdot 5\cdot 7\cdot 127}\func{mod}1$.
\end{quote}

\noindent Since the coboundary $W$ is spin we can use the Spin
characteristic classes to replace the Pontryagin classes by the transition
(9.4) to get the following simpler expressions.

\bigskip

\noindent \textbf{Proposition 10.3.} \textsl{In term of the Spin
characteristic classes, the Eells and Kuiper invariants }$\mu _{k}(M)$%
\textsl{,} $k=2,3$ \textsl{or} $4,$ \textsl{have the following expressions}

\begin{quote}
$\mu _{2}(M)\equiv \frac{(q_{1}^{2}-\sigma )[W]}{2^{5}\cdot 7}\func{mod}1$

$\mu _{3}(M)\equiv \frac{(2(q_{1}q_{2}-q_{1}^{3})-3\sigma )[W]}{2^{7}\cdot
(2^{5}-1)\cdot 3}\func{mod}1$

$\mu _{4}(M)\equiv \frac{%
(6q_{1}q_{3}+5q_{2}^{2}-40q_{1}^{2}q_{2}+17q_{1}^{4}-45\sigma )[W]}{%
2^{9}(2^{7}-1)}\func{mod}1$.$\square $
\end{quote}

\textbf{10.4. The Rohlin type formula.} For an $4m$ dimensional oriented
manifold $M$ denote by $\sigma _{M}$ the signature of the intersection form $%
I_{M}:H^{2m}(M)\rightarrow \mathbb{Z}$, $I_{M}(x)=\left\langle
x^{2},[M]\right\rangle $, on the middle dimensional integral cohomology.

By the year 1950 there had no example of topological manifolds that are not
smoothable. Nevertheless, Whitehead \cite[1949]{Wh} had constructed for each
integral unimodular symmetric matrix $A$ a simply--connected $4$ dimensional
topological manifold $M$, whose the intersection form $I_{M}$ is precisely
given by $A$. In contrast, the following result of Rokhlin \cite{Ro,Fr}
singles out a severe gap between the topological and smooth categories,
detectable by the elementary and topological invariant $\sigma _{M}$.

\bigskip

\noindent \textbf{Theorem 10.4 (Rokhlin, 1952).} \textsl{The signature }$%
\sigma _{M}$ \textsl{of a }$4$\textsl{--dimensional smooth spin manifold }$M$
\textsl{must be divisible by }$16$\textsl{.}$\square $

\bigskip

To extend Rokhlin's result to the higher dimensional settings we resort to
the $\widehat{A}$ genus $\alpha _{m}$ and the $L$ genus $\tau _{m}$ of $4m$
dimensional smooth oriented manifolds $M^{4m}$, which are certain
polynomials in the Pontryagin classes $p_{1},\cdots ,p_{m}$ of $M$ with
homogeneous degree $4m$

\begin{quote}
$\alpha _{m}=a_{m}\cdot p_{m}+l_{m}(p_{1},\cdots ,p_{m-1})$ and

$\tau _{m}=b_{m}\cdot p_{m}+k_{m}(p_{1},\cdots ,p_{m-1})$ \cite[\S 19]{MS},
\end{quote}

\noindent where $a_{m}$ and $b_{m}$ are certain non--zero rationals, and
where $l_{m}$ and $k_{m}$ are certain polynomials in $p_{1},\cdots ,p_{m-1}$
with rational coefficients. These allow us to eliminate the top degree
Pontryagin class $p_{m}$ to obtain the following expression of $\tau _{m}$
without involving $p_{m}$:

\begin{enumerate}
\item[(10.2)] $\tau _{m}=$ $\frac{b_{m}}{a_{m}}(\alpha
_{m}-l_{m}(p_{1},\cdots ,p_{m-1}))+k_{m}(p_{1},\cdots ,p_{m-1})$.
\end{enumerate}

\noindent For the geometric implications of the polynomials $\alpha _{m}$
and $\tau _{m}$ we recall the following classical results.

\bigskip

\noindent \textbf{Theorem 10.5 (Hirzebruch }\cite{BH}\textbf{) }\textsl{The }%
$L$\textsl{--genus }$\tau _{m}$\textsl{\ of }$M^{4m}$\textsl{\ equals to the
signature }$\sigma _{M}$ \textsl{of the intersection form on }$%
H^{2m}(M^{4m}) $\textsl{.}$\square $

\bigskip

\noindent \textbf{Theorem 10.6 (Borel--Hirzebruch \cite{BH})} \textsl{The }$%
\widehat{A}$ \textsl{genus} $\alpha _{m}$\textsl{\ of a spin manifold }$%
M^{4m}$ \textsl{is an integer, and is an even integer when }$m$\textsl{\ is
odd.}$\square $

\bigskip

\noindent \textbf{Theorem 10.7 (Gromov, Lawson and Stolz \cite{L,S0}).} 
\textsl{If }$M^{4m}$\textsl{\ is a simply connected spin manifold with }$m>1$%
\textsl{, then }$M^{4m}$\textsl{\ admits a metric with positive scalar
curvature if and only if }$\alpha _{m}=0$\textsl{.}$\square $

\bigskip

Precisely, for $1\leq m\leq 4$ the polynomials $\alpha _{m}$ and $\tau _{m}$
are, respectively,

\begin{quote}
$\alpha _{1}=-\frac{1}{24}p_{1}$;

$\alpha _{2}=\frac{1}{2^{7}\cdot 3^{2}\cdot 5}(-4p_{2}+7p_{1}^{2});$

$\alpha _{3}=\frac{1}{2^{10}\cdot 3^{3}\cdot 5\cdot 7}%
(-16p_{3}+44p_{2}p_{1}-31p_{1}^{3});$

$\alpha _{4}=\frac{1}{2^{15}\cdot 5^{2}\cdot 3^{4}\cdot 7}%
(-192p_{4}+512\cdot p_{1}p_{3}+208p_{2}^{2}-904p_{1}^{2}p_{2}+381p_{1}^{4}),$
\end{quote}

\noindent and

\begin{quote}
$\tau _{1}=\frac{1}{3}p_{1}$;

$\tau _{2}=\frac{1}{3^{2}\cdot 5}(7p_{2}-p_{1}^{2})$;

$\tau _{3}=\frac{1}{3^{3}\cdot 5\cdot 7}(62p_{3}-13p_{2}p_{1}+2p_{1}^{3})$;

$\tau _{4}=\frac{1}{3^{4}\cdot 5^{2}\cdot 7}(381p_{4}-71\cdot
p_{1}p_{3}-19p_{2}^{2}+22p_{1}^{2}p_{2}-3p_{1}^{4})$.
\end{quote}

\noindent Assume now that our manifold $M$ is spin (i.e. $w_{2}(M)=0$). Then
the formulae in (9.4) is applicable to replace the Pontryagin classes $%
p_{1},\cdots ,p_{m-1}$ in (10.2) to yield the following simpler formulae of
the signature $\sigma _{M}=\tau _{m}$ in the Spin characteristic classes.

\bigskip

\noindent \textbf{Theorem 10.8. }\textsl{In accordance to }$m=1,2,3$\textsl{%
\ and }$4$\textsl{,\ the signature }$\sigma _{M}$\textsl{\ of a smooth spin
manifold }$M^{4m}$ \textsl{is given, respectively, by}

\begin{quote}
$\sigma _{M}=-2^{3}\cdot \alpha _{1}$;

$\sigma _{M}=q_{1}^{2}-2^{5}\cdot (2^{3}-1)\cdot \alpha _{2}$;

$\sigma _{M}=\frac{2}{3}(q_{1}q_{2}-q_{1}^{3})-2^{7}\cdot (2^{5}-1)\cdot
\alpha _{3}$;

$\sigma _{M}=\frac{2}{3\cdot 5}q_{1}q_{3}+\frac{1}{3^{2}}q_{2}^{2}-\frac{%
2^{3}}{3^{2}}q_{1}^{2}q_{2}+\frac{17}{3^{2}\cdot 5}q_{1}^{4}-2^{9}(2^{7}-1)%
\cdot \alpha _{4}$,
\end{quote}

\noindent \textsl{where }$\alpha _{1}\equiv \alpha _{3}\equiv 0\func{mod}2$%
\textsl{.}$\square $

\bigskip

\noindent \textbf{Example 10.9.} For a smooth spin manifold $M$ one has by
Theorems 10.5, 10.6 and 10.8 that

\begin{quote}
i) if $\dim M=4$, then $\sigma _{M}\equiv 0\func{mod}2^{4}$;

ii) if $\dim M=8$, then $\sigma _{M}\equiv q_{1}^{2}\func{mod}2^{5}\cdot
(2^{3}-1)$\textsl{.}

iii) if $\dim M=12$, then $\sigma _{M}\equiv \frac{2}{3}%
(q_{1}q_{2}-q_{1}^{3})\func{mod}2^{8}\cdot (2^{5}-1)$\textsl{,}
\end{quote}

\noindent where assertion i) is identical to Theorem 10.4. For this reason
we may call the formulae in Theorem 10.8 \textsl{the Rokhlin type formulae}
of spin manifolds.

Since the string group $String(n)$ is the $3$--connected cover of $Spin(n)$,
a spin manifold is \textsl{string} \cite{S} if and only if its first spin
characteristic class $q_{1}$ vanishes. In this case the second spin
characteristic class $q_{2}$ has been shown to be divisible by $3$ \cite{LD}%
. Thus, for a smooth string manifold $M$ we have by Theorem 10.8 the
following Rokhlin type formulae.

\begin{quote}
a) if $\dim M=8$, then $\sigma _{M}\equiv 0\func{mod}2^{5}\cdot (2^{3}-1)$%
\textsl{.}

b) if $\dim M=12$, then $\sigma _{M}\equiv 0\func{mod}2^{8}\cdot (2^{5}-1)$%
\textsl{,}

c) if $\dim M=16$, then $\sigma _{M}\equiv (\frac{1}{3}q_{2})^{2}\func{mod}%
2^{9}\cdot (2^{7}-1)$.
\end{quote}

\noindent In addition, by Theorem 10.7 if $M$\ is a simply connected, then $%
M $\ admits a metric with positive scalar curvature if and only if

\begin{quote}
$\sigma _{M}=0,0$ or $(\frac{1}{3}q_{2})^{2}$ in accordance to $\dim M=8,12$
or $16$.$\square $
\end{quote}

\textbf{10.5. The existence of smooth structure on triangulable manifolds.}
Without involving the top degree Pontryagin class $p_{m}$ the Rokhlin type
formulae in Theorem 10.8 is ready to apply to study of the existence problem
of smooth structures on certain $4m$ dimensional triangulable manifolds. To
provide such examples in the case $m=2$ we need the following notation.

\bigskip

\noindent \textbf{Definition 10.10.} For a unimodular symmetric integral
matrix $A=(a_{ij})_{n\times n}$ of rank $n$, and a sequence $b=(b_{1},\cdots
,b_{n})$ of integers with length $n$, the pair $(A,b)$ is called a\textsl{\
Wall pair} if the following congruences are satisfied

\begin{enumerate}
\item[(10.3)] $a_{ii}\equiv b_{i}\func{mod}2,1\leq i\leq n$.$\square $
\end{enumerate}

\bigskip

Let $D^{8}$ be the unit disk on the $8$ dimensional Euclidean space $\mathbb{%
R}^{8}$. The following result is due to C.T.C. Wall \cite{Wa}.

\bigskip

\noindent \textbf{Theorem 10.11. }\textsl{For each Wall pair }$(A,b)$\textsl{%
\ with }$A=(a_{ij})_{n\times n}$ \textsl{and }$b=(b_{1},\cdots ,b_{n})$%
\textsl{,} \textsl{there exists a closed }$8$\textsl{\ dimensional} \textsl{%
topological manifold }$M$\textsl{\ that satisfies the following properties}

\textsl{i) }$M$\textsl{\ admits a decomposition }$M=W\cup _{h}D^{8}$\textsl{%
, where }$W$ \textsl{is a }$3$\ \textsl{connected} \textsl{smooth manifold
with boundary }$\partial W$ \textsl{a homotopy }$7$\textsl{--sphere, and
where }$h:\partial W\rightarrow \partial D^{8}$ \textsl{is a homeomorphism;}

\textsl{ii) there is a basis }$\left\{ x_{1},\cdots ,x_{n}\right\} $\textsl{%
\ on }$H^{4}(M)$\textsl{\ so that }$x_{i}\cup x_{j}=a_{i,j}\cdot \omega _{M}$%
\textsl{, where }$\omega _{M}\in H^{8}(M)$ \textsl{is an orientation class on%
} $M$\textsl{;}

\textsl{iii) the first Spin characteristic class }$q_{1}$\textsl{\ of }$M$%
\textsl{\ is well defined (by i)), and is determined by }$b$\textsl{\ as }

\begin{quote}
$\qquad q_{1}=b_{1}x_{1}+\cdots +b_{n}x_{n}\in H^{4}(M)$\textsl{.}
\end{quote}

\textsl{Furthermore, if }$(A^{\prime },b^{\prime })$\textsl{\ is a second
Wall pair, then the associated manifold }$M^{\prime }$\textsl{\ is
combinatorially homeomorphic to }$M$\textsl{\ (in the sense of \cite{Wa}) if
and only if there exists an integer matrix }$P=(p_{i,j})_{n\times n}$\textsl{%
\ so that }

\begin{quote}
$P^{\tau }AP=A^{\prime }$ \textsl{and} $bP=b^{\prime }$\textsl{, }
\end{quote}

\noindent \textsl{where }$P^{\tau }$\textsl{\ denotes the transpose of the
matrix }$P$\textsl{.}$\square $

\bigskip

\noindent \textbf{Remark 10.12.} In \cite{Wa} Wall classified the
combinatorial homeomorphism types of all the $(n-1)$ connected $2n$
dimensional manifolds $M$ that are smooth off one point $o\in M$. The result
in Theorem 10.11 corresponds to the case $n=4$.

It is known that for a $4$--dimensional real vector bundle $\xi $ on the $4$
dimensional sphere $S^{4}$ the difference $2e(\xi )-p_{1}(\xi )$ (resp. $%
e(\xi )-q_{1}(\xi )$) is divisible by $4$ (resp. by $2$), where $e(\xi )$ is
the Euler class of $\xi $ \cite[Lemma 20.10]{MS}. In Theorem 10.11 the
necessity of the Wall condition (10.3) is governed by the following
geometric fact. According to Haefliger \cite{Ha}, for the manifold $W$ in i)
of Theorem 10.11, there exist $n$ smooth embeddings

\begin{quote}
$\iota _{i}:S^{4}\rightarrow W$, $1\leq i\leq n$,
\end{quote}

\noindent so that the Kronecker dual of the cycle classes $\iota _{i\ast
}[S^{4}]\in H_{4}(M)$ is the basis $\left\{ x_{1},\cdots ,x_{n}\right\} $%
\textsl{\ }on\textsl{\ }$H^{4}(M)$. Then, the matrix $A$ is the intersection
form on $H^{4}(M)$ corresponding to the basis, while the normal bundle $%
\gamma _{i}$ of the embedding $\iota _{i}$ is related to the pair $(A,b)$ by
the relations

\begin{quote}
$(e(\gamma _{i}),q_{1}(\gamma _{i}))=(a_{ii}\cdot \omega ,b_{i}\cdot \omega
) $, $1\leq i\leq n,$
\end{quote}

\noindent where $\omega $ is the orientation class on $S^{4}$ that
corresponds $x_{i}$ via $\iota _{i}$.

For an $8$--dimensional manifold $M$ associated to a Wall pair $(A,b)$
properties ii) and iii) of Theorem 10.9 imply, respectively, that

\begin{enumerate}
\item[(10.4)] $\sigma _{M}=sign(A)$ and $q_{1}(M)^{2}=bAb^{\tau }$,
\end{enumerate}

\noindent where $b^{\tau }$\textsl{\ }denotes the transpose of the row
vector $b$.$\square $

\bigskip

Concerning the manifold $M$ associated to a Wall pair $(A,b)$ a natural
question is whether there exists a smooth structure that extends the given
one on $W$. For the special case $A=(1)_{1\times 1}$ this question has been
studied by Milnor \cite{M}, Eells and Kuiper \cite[\S 6]{EK} in their
calculation on the group $\Theta _{7}$ of homotopy $7$ spheres. We extend
their calculations in the following results.

\bigskip

\noindent \textbf{Theorem 10.13. }\textsl{Let }$M^{8}$\textsl{\ be the
manifold associated to a Wall pair }$(A,b)$\textsl{. There exists a smooth
structure} \textsl{on }$M^{8}$\textsl{\ extending the one on }$W$\textsl{\
if and only if}

\begin{enumerate}
\item[(10.5)] $sign(A)\equiv bAb^{\tau }\func{mod}2^{5}\cdot (2^{3}-1)$%
\textsl{.}
\end{enumerate}

\bigskip

\noindent \textbf{Proof.} The necessity of (10.5) comes from ii) of Example
10.9. The sufficiency is verified by computing with the Eells--Kuiper $\mu $
invariant \cite[formula (11)]{EK} of the boundary $\partial W$, which by
Proposition 10.3 reads

\begin{quote}
$\mu (\partial W)\equiv \frac{4bAb^{\tau }-4sign(A)}{2^{7}\cdot (2^{3}-1)}%
\equiv \frac{bAb^{\tau }-sign(A)}{2^{5}\cdot (2^{3}-1)}\func{mod}1$.$\square 
$
\end{quote}

Theorem 10.13 has several direct, but notable consequences. A theorem of
Kervaire states that there exist a $10$ dimensional manifold which do not
admit any smooth structure \cite{K}. Eells and Kuiper provided further
examples which have the same cohomology ring as that of the projective plane 
\cite{EK1}. Theorem 10.13 implies that

\bigskip

\noindent \textbf{Corollary 10.14. }\textsl{If }$(A,b)$\textsl{\ is a Wall
pair so that }

\begin{quote}
$sign(A)\neq bAb^{\tau }\func{mod}2^{5}\cdot (2^{3}-1)$\textsl{, }
\end{quote}

\noindent \textsl{then the corresponding manifold} $M$ \textsl{does not
admit any smooth structure.}$\square $

\bigskip

Conversely, for those $M$ which do admit smooth structures, their total Spin
characteristic class $q(M)$ can be determined completely.

\bigskip

\noindent \textbf{Corollary} \textbf{10.15.} \textsl{If the manifold} $M$ 
\textsl{associated to a Wall pair }$(A,b)$ \textsl{admits a smooth
structure, then its total Spin characteristic class is}

\begin{enumerate}
\item[(10.6)] $q(M)=1+(b_{1}x_{1}+\cdots +b_{1}x_{n})+\frac{3(15\cdot
sign(A)-bAb^{\tau })}{2\cdot (2^{3}-1)}\cdot \omega _{M}$\textsl{,}
\end{enumerate}

\noindent \textsl{where }$2\cdot (2^{3}-1)$\textsl{\ divides }$15\cdot
sign(A)-bAb^{\tau }$\textsl{\ by ii) of Example 10.9..}

\bigskip

\noindent \textbf{Proof.} With $\tau _{2}=\sigma _{M}=sign(A)$ and $%
p_{1}^{2}=4bAb^{\tau }$ the formula $\tau _{2}=\frac{7p_{2}-p_{1}^{2}}{%
3^{2}\cdot 5}$ implies that

\begin{quote}
$p_{2}=\frac{45\cdot sign(A)+4bAb^{\tau }}{7}\cdot \omega _{M}$.
\end{quote}

\noindent From $p_{2}=2q_{2}+q_{1}^{2}$ by (9.4) we get $q_{2}=\frac{%
3(15\cdot sign(A)-bAb^{\tau })}{2\cdot (2^{3}-1)}\cdot \omega _{M}$.$\square 
$

\bigskip

In view of the formula of $q_{2}$ in (10.6) Theorem 10.7 implies that

\bigskip

\noindent \textbf{Corollary 10.16. }\textsl{For} \textsl{a manifold} $M$ 
\textsl{associated to a Wall pair }$(A,b)$\textsl{\ the following statements
are equivalent:}

\begin{quote}
\textsl{i) }$M$\textsl{\ is smoothable and has a metric with positive scalar
curvature;}

\textsl{ii)} $sign(A)=bAb^{\tau }$\textsl{;}

\textsl{iii)} $q_{2}=3sign(A)\cdot \omega _{M}$\textsl{.}$\square $
\end{quote}

\bigskip

\noindent \textbf{Example 10.17. }Corollary 10.16 reduces the problem of
finding all the $3$--connected and $8$-- dimensional smooth manifolds that
have a metric with positive scalar curvature to the arithmetic problem of
finding those Wall pairs $(A,b)$ satisfying the quadratic equation

\begin{enumerate}
\item[(10.7)] $sign(A)=bAb^{\tau }$.
\end{enumerate}

Consider a Wall pair $(A,b)$ with $A$ the identity matrix $I_{n}$ of rank $n$%
, and with $b=(2k_{1}+1,\cdots ,2k_{n}+1)$, $k_{i}\in \mathbb{Z}$ (see
(10.3)). The equation (10.7) is

\begin{quote}
$n=(2k_{1}+1)^{2}+\cdots +(2k_{n}+1)^{2}$.
\end{quote}

\noindent It implies that the corresponding $M$\textsl{\ }is smooth and
admits a metric with positive scalar curvature, if and only if $M$ is
combinatorially homeomorphic \cite{Wa}\textsl{\ }to $\mathbb{H}P^{2}\#\cdots
\#\mathbb{H}P^{2}$, the connected sum of $n$ copies of the projective plan $%
\mathbb{H}P^{2}$.

Consider next a Wall pair $(A,b)$ with $A=$ $\left( 
\begin{array}{cc}
0 & 1 \\ 
1 & 0%
\end{array}%
\right) $ and $b=(2k_{1},2k_{2})$, $k_{i}\in \mathbb{Z}$ (see (10.3)). The
equation (10.7) turns to be $k_{1}\cdot k_{2}=0$. It implies that the
corresponding manifold $M$\textsl{\ }is smooth and admits a metric with
positive scalar curvature, if and only if $M$ is combinatorially
homeomorphic to the spherical bundle $S(\xi )$ of a $5$ dimensional
Euclidean bundle $\xi $ over $S^{4}$. Such manifolds $S(\xi )$ are
classified by the homotopy group $\pi _{3}(SO(5))$.

Concerning this topic general cases will be studied in the sequel work \cite%
{DL}.$\square $

\bigskip

\textbf{Acknowledgement.} The author would like to thank Fei Han, Ruizhi
Huang for discussion on the integral lifts of Wu--classes \cite{HS,DHH} (see
Section \S 9.2), and to Yang Su for informing him the recent questions in
mathoverflow \cite{Web,Web1} about the integral cohomology of the
classifying space $B_{Spin(n)}$.

\end{document}